\documentclass[12pt]{amsart}
\usepackage{amssymb,amsmath, amsthm,latexsym}
\usepackage{a4wide}

\theoremstyle{plain}

\newtheorem{lem}{Lemma}[section]
\newtheorem{theo}[lem]{Theorem}
\newtheorem{prop}[lem]{Proposition}
\newtheorem{corollary}[lem]{Corollary}
\newtheorem{remark}[lem]{Remark}

\parskip=\medskipamount
\font\k=cmr7
.7

  \newcommand {\di}{\mbox{\k disc}}
  \newcommand {\reg}{\mbox{\k reg}}
  
  \newcommand {\C}{{\mathbb C}}
  \newcommand {\N}{{\mathbb N}}
  \newcommand {\R}{{\mathbb R}}
  \newcommand {\Z}{{\mathbb Z}}
  \newcommand {\Q}{{\mathbb Q}}
  \newcommand {\A}{{\mathbb A}}
  \newcommand {\af}{{\mathfrak a}}

  \newcommand {\ho}{{\mathfrak o}}
  \newcommand {\qf}{{\mathfrak q}}
  \newcommand {\pg}{{\mathfrak p}}
  
  \newcommand {\Pg}{{\mathfrak P}}
 \newcommand {\mX}{{\mathfrak X}}
 \newcommand {\mM}{{\mathfrak M}}
 
 \newcommand {\mO}{{\mathfrak O}}
  
\renewcommand {\H}{{\mathcal H}}
  
  \newcommand {\Co}{{\mathcal C}}
  \newcommand {\cO}{{\mathcal O}}
  
\renewcommand {\L}{{\mathcal L}}

 \newcommand {\cP}{{\mathcal P}}
 \newcommand {\cF}{{\mathcal F}}
 \newcommand {\cL}{{\mathcal L}}
 \newcommand {\cA}{{\mathcal A}}
\newcommand {\ba}{\backslash}
 \newcommand {\ov}{\overline}
  \newcommand {\bs}{{\bf s}}
 \newcommand {\bu}{{\bf u}}
 
 \newcommand {\bk}{{\bf k}}
\newcommand{\norm}[1]{\lVert#1\rVert}
\newcommand{\abs}[1]{\left|{#1}\right|}

\renewcommand{\Re}{\operatorname{Re}}

\newcommand{\sign}{\operatorname{sign}}
\newcommand{\tr}{\operatorname{tr}}

\newcommand{\Hom}{\operatorname{Hom}}
\newcommand{\Ind}{\operatorname{Ind}}
\newcommand{\I}{\operatorname{I}}
\newcommand{\vol}{\operatorname{vol}}
\newcommand{\Res}{\operatorname{Res}}

\newcommand{\GL}{\operatorname{GL}}
\newcommand{\id}{\operatorname{id}}
\newcommand{\rO}{\operatorname{O}}

\newcommand{\U}{\operatorname{U}}

\begin{document}
\title[Arthur Trace formula]{\Large\bf On the absolute convergence of the  spectral
 side
 of the Arthur trace formula for $\GL_n$}
\date{\today}
\footnotetext{$^\dagger$partially supported by NSF grant DMS-0070561}

\address{Universit\"at Bonn\\
Mathematisches Institut\\
Beringstrasse 1\\
D -- 53115 Bonn, Germany}
\email{mueller@math.uni-bonn.de}
\address{Department of Mathematics\\
310 Malott Hall\\
Cornell University\\
Ithaca, NY 14853-4201\\
U.S.A.}
\email{speh@math.cornell.edu}
\address{Einstein Institute of Mathematics\\
 The Hebrew University of Jerusa\-lem\\
Jerusalem 91904\\
Israel}
\email{erezla@math.huji.ac.il}

\maketitle

\centerline{\sc Werner M\"uller and Birgit Speh$^\dagger$}
\vskip2pt
\centerline{\sc with Appendix by}
\centerline{\sc Erez M. Lapid}

\begin{abstract}
Let $G$ be a reductive algebraic group defined over $\Q$ and let $\A$ be the
 ring of ad\`eles of $\Q$. The spectral side
of the Arthur trace formula for $G$ is a sum of distributions on $G(\A)^1$
which are defined in terms of truncated Eisenstein series. In general, the
spectral side is only known to be  conditionally convergent. In this paper we
 prove
that for $\GL_n$,  the spectral side of the
 trace formula is absolutely convergent. 
\end{abstract}

\setcounter{section}{-1}
\section{Introduction}
Let $E$ be a number field and let $G$ be a connected reductive algebraic
 group over $E$.
 Let $\A$ be the ring of ad\`eles of $E$
and let $G(\A)$ be the group of points of $G$ with values in $\A$. 
Let $G(\A)^1$ be the intersection of the kernels of the
maps $x\mapsto|\xi(x)|$, $x\in G(\A)$, where $\xi$ ranges over the group
$X(G)_E$ of characters of $G$ defined over $E$. Then the (noninvariant)
trace formula of Arthur is an identity
$$\sum_{\ho\in\mO}J_\ho(f)=\sum_{\chi\in\mX}J_\chi(f),\quad f\in C_c^\infty(G(\A)^1),$$
between distributions on $G(\A)^1$. The left hand side is  the
{\it geometric side} and the right hand side the {\it spectral side} of the 
trace formula. 

In this paper we are concerned with the spectral side of the
trace formula. The distributions $J_\chi$
are initially defined in terms of truncated Eisenstein series. They are 
parametrized 
by the set  of cuspidal data $\mX$ which consists of the Weyl group orbits of
pairs $(M_B,r_B)$, where $M_B$ is the Levi component of a standard parabolic
subgroup and $r_B$ is an irreducible cuspidal automorphic representation 
of $M_B(\A)^1$. In the fine $\chi$-expansion of the spectral side the inner
products of truncated Eisenstein series are replaced by terms containing 
generalized logarithmic derivatives of intertwining operators. This leads to
an integral-series that is only known to be conditionally convergent. It is
an open problem to prove that the fine $\chi$-expansion is absolutely
convergent and the main purpose of this paper is to
settle this problem for the group $\GL_n$. 

To explain our results in more detail, we need to introduce some notation.
We fix a Levi component $M_0$  of a minimal parabolic subgroup  $P_0$  of
$G$. We assume that all parabolic 
subgroups considered in this paper contain $M_0$.
Let $P$ be a parabolic subgroup 
of $G$, defined over $E$, with unipotent radical $N_P$. Let $M_P$ be the 
unique Levi component of $P$ which contains $M_0$. We denote the split 
component of the center of $M_P$ by $A_P$ and its Lie algebra by $\af_P$. 
For parabolic groups $P\subset Q$ there is a natural surjective map $\af_P\to
\af_Q$ whose kernel we will denote by $\af_P^Q$.
Let $\cA^2(P)$ be the space of  automorphic forms on
$N_P(\A)M_P(E)\backslash G(\A)$ which are square-integrable modulo 
$A_{P,\Q}(\R)^0$, where $A_{P,\Q}$ is the split component of the center of the
group obtained from $M_P$ by restricting scalars from $E$ to $\Q$.
Let $Q$ be another parabolic subgroup of $G$, defined over $E$, with Levi
component $M_Q$, split component $A_Q$ and corresponding Lie algebra $\af_Q$.
Let $W(\af_P,\af_Q)$ be the set of all linear isomorphisms from $\af_P$ to 
$\af_Q$
which are restrictions of elements of the Weyl group $W(A_0)$.  The theory of 
Eisenstein series associates to each $s\in W(\af_P,\af_Q)$ an intertwining 
operator
$$M_{Q|P}(s,\lambda): \cA^2(P)\rightarrow \cA^2(Q),\quad \lambda\in\af_{P,\C}^*,$$
which for $\Re(\lambda)$ in a certain chamber, can be defined by an 
absolutely convergent integral and admits an analytic continuation to a 
meromorphic function of $\lambda\in\af^*_{P,\C}$. Set
$$M_{Q|P}(\lambda):=M_{Q|P}(1,\lambda).$$
Let $\Pi(M_P(\A)^1)$ be the
set of equivalence classes of irreducible unitary representations of 
$M_P(\A)^1$. Let $\chi\in\mX$ and 
$\pi\in\Pi(M_P(\A)^1)$. Then $(\chi,\pi)$ singles out a certain subspace
$\cA^2_{\chi,\pi}(P)$ of $\cA^2(P)$ \cite[p.1249]{A3}.
 Let $\ov \cA^2_{\chi,\pi}(P)$ be the Hilbert 
space completion of $\cA^2_{\chi,\pi}(P)$ with respect to the canonical inner
product. For each $\lambda\in\af_{P,\C}^*$ we have an induced representation
 $\rho_{\chi,\pi}(P,\lambda)$  of $G(\A)$ in $\ov \cA^2_{\chi,\pi}(P)$.

For each
Levi subgroup $L$ let $\cP(L)$ be the set of all parabolic subgroups with
Levi component $L$. If $P$ is a parabolic subgroup, let $\Delta_P$ denote the
set of simple roots of $(P,A_P)$. Let $L$ be a Levi subgroup which contains 
$M_P$. Set
\begin{equation*}
\begin{split}
\mM_L&(P,\lambda)=\\
&\lim_{\Lambda\to0}\left(\sum_{Q_1\in\cP(L)}
\vol(\af_{Q_1}^G/\Z(\Delta^\vee_{Q_1}))M_{Q|P}(\lambda)^{-1}
\frac{M_{Q|P}(\lambda+\Lambda)}{\prod_{\alpha\in\Delta_{Q_1}}\Lambda(\alpha^\vee)}\right),
\end{split}
\end{equation*}
where $\lambda$ and $\Lambda$ are constrained to lie in $i\af_L^*$, and for 
each $Q_1\in\cP(L)$, $Q$ is a group in $\cP(M_P)$ which is contained in $Q_1$.
Then $\mM_L(P,\lambda)$ is an unbounded operator which acts
on the Hilbert space $\ov \cA^2_{\chi,\pi}(P)$.  
In the special case that $L=M$ 
and $\dim\af_L^G=1$, the operator $\mM_L(P,\lambda)$ has a simple description. 
 Let $P$ be a parabolic subgroup with Levi component $M$. 
 Let $\alpha$ be the unique simple root of  $(P,A_P)$ and let $\tilde\omega$
be the element in $(\af_M^G)^*$ such that $\tilde\omega(\alpha^\vee)=1$. Let
$\ov P$ be the opposite parabolic group of $P$.  Then
$$\mM_L(P,z\tilde\omega)=-\vol(\af_M^G/\Z\alpha^\vee)
M_{\ov P|P}(z\tilde\omega)^{-1}\cdot
\frac{d}{dz}M_{\ov P|P}(z\tilde\omega).$$

Let $f\in C_c^\infty(G(\A)^1)$. Then Arthur
\cite[Theorem 8.2]{A4} proved that $J_\chi(f)$ equals the sum over Levi 
subgroups $M$ containing $M_0$, over $L$ containing $M$, over $\pi\in
\Pi(M(\A)^1)$, and over $s\in W^L(\af_M)_{\reg}$, a certain subset of the 
Weyl group, of the product of
\begin{equation*}
|W_0^M||W_0|^{-1}|\det(s-1)_{\af^L_M}|^{-1}|\cP(M)|^{-1}
\end{equation*}
a  factor to which we need not pay too much attention , and of
\begin{equation}\label{0.1}
\int_{i\af_L^*/i\af_G^*}\sum_{P\in\cP(M)}\tr(\mM_L(P,\lambda)M_{P|P}(s,0)
\rho_{\chi,\pi}(P,\lambda,f))\;d\lambda.
\end{equation}
So far, it is only known that $\sum_{\chi\in\mX}|J_\chi(f)|<\infty$
and the goal is to show that the integral--sum obtained by summing 
(\ref{0.1}) over $\chi\in\mX$ and $\pi\in\Pi(M(\A)^1)$ is absolutely 
convergent with respect to the trace norm.
For a given Levi subgroup $M$ let $\L(M)$
be the set of all Levi subgroups $L$ with $M\subset L$. Put 
$M(P,s)=M_{P|P}(s,0)$. Denote by $\parallel T\parallel_1$ the trace norm
of a trace class operator $T$.  Let $\Co^1(G(\A)^1)$ be the space of 
integrable rapidly decreasing functions on $G(\A)^1$ (see  \cite[\S1.3]{Mu4}
 for its definition).
Then our main result is the following theorem.

\begin{theo}\label{th0.1} 
Let $G=\GL_n$. 
Then the sum over all $M\in{\L}(M_0)$, $L\in{\L}(M)$,
 $\chi\in\mX$,  $\pi\in\Pi(M(\A)^1)$,  and $s\in W^L(\af_M)_{\reg}$ 
of the product of 
\begin{equation*}
|W^M_0| |W_0|^{-1} |\det(s-1)_{\af^L_M}|^{-1} 
\end{equation*}
with
\begin{equation*}
\int_{i\af^*_L/i \af^*_G}  |\cP(M)|^{-1}\sum_{P\in\cP(M)}
\parallel \mM_L
(P,\lambda)M(P,s)\rho_{\chi,\pi}(P,\lambda,f)\parallel_1d\lambda 
\end{equation*}
is convergent for all $f\in\Co^1(G(\A)^1)$.
\end{theo}

By Theorem \ref{th0.1}, the spectral side for $\GL_n$ 
can now be rewritten in the 
following way. Denote by $\Pi_{\di}(M(\A)^1)$ the set of all $\pi\in
\Pi(M(\A)^1)$ which are equivalent to an irreducible subrepresentation of the
regular representation of $M(\A)^1$ in $L^2(M(E)\ba M(\A)^1)$. As in 
Section 7 of \cite{A3}, we shall identify any representation of $M(\A)^1$
with a representation of $M(\A)$ which is trivial on 
$A_{M,\Q}(\R)^0$, where $A_{M,\Q}$ is the split component of the center of 
the group $\Res_{E/\Q}\GL_n$ obtained from $G$ by restricting scalars from 
$E$ to $\Q$. For any parabolic group $P$, let 
$\cA^2_\pi(P)=\oplus_\chi\cA^2_{\chi,\pi}(P)$
 and for $\lambda\in\af_{P,\C}^*$, let 
$\rho_\pi(P,\lambda)$ be
the induced representation of $G(\A)$ in $\overline\cA^2_\pi(P)$, the Hilbert 
space completion of $\cA^2_\pi(P)$.
Given $M\in\cL$, $L\in\cL(M)$, $P\in\cP(M)$,  $s\in W^L(\af_M)_{\reg}$ and a
function $f\in\Co^1(G(\A)^1)$, let
\begin{equation*}
\begin{split}
J^L_{M,P}&(f,s)\\
&=\sum_{\pi\in\Pi_{\di}(M(\A)^1)}
\int_{i\af_L^*/i\af_G^*}\tr(\mM_L(P,\lambda)M_{P|P}(s,0)
\rho_\pi(P,\lambda,f))\;d\lambda.
\end{split}
\end{equation*}

By Theorem \ref{th0.1} this integral-series is absolutely convergent with 
respect to the trace norm.
 Furthermore for $M\in\cL$ and $s\in W^L(\af_M)_{\reg}$ set
$$a_{M,s}=|\cP(M)|^{-1}|W^M_0||W_0|^{-1}|\det(s-1)_{\af^L_M}|^{-1}.$$
Then for all functions $f$ in $\Co^1(G(\A)^1)$, the spectral  side of
the Arthur trace formula equals
$$\sum_{M\in\cL}\sum_{L\in\cL(M)}\sum_{P\in\cP(M)}\sum_{s\in W^L(\af_M)_{\reg}}
 a_{M,s} J^L_{M,P}(f,s).$$
Note that all sums in this expression are finite.

We shall now explain the main steps of the proof of Theorem \ref{0.1}. The
 proof relies on Theorem 0.1 of \cite{Mu4}. 
In this theorem  the absolute convergence of the 
spectral side of the trace formula has been reduced to a problem about local
components of automorphic representations. So the main issue of the present
paper is to verify that  for $\GL_n$, the assumptions of Theorem 0.1 of 
\cite{Mu4} are satisfied.

Let $M=\GL_{n_1}\times\cdots\times\GL_{n_r}$ be a standard Levi subgroup and
let $P,Q\in\cP(M)$. We shall identify $\af_M^*$ with $\R^r$. Given a place
$v$ of $E$ and an irreducible unitary representation $\pi_v=\otimes_{i=1}^r
\pi_{v,i}$ of $M(E_v)$, let $J_{Q|P}(\pi_v,\bs)$, $\bs\in\C^r$,  be the
 local intertwining operator between the induced representations 
$I^G_P(\pi_{v}[\bs])$ and $I^G_Q(\pi_{v}[\bs])$, where $\bs=(s_1,...,s_r)$
and $\pi_v[\bs]=\otimes_i(\pi_{v,i}|\det|^{s_i})$. It follows from results
of Shahidi \cite{Sh5} that there exist normalizing factors 
$r_{Q|P}(\pi_v,\bs)$, which are
defined in terms of Rankin-Selberg $L$-functions,  such that the normalized
 intertwining operators
$$R_{Q|P}(\pi_v,\bs)=r_{Q|P}(\pi_v,\bs)^{-1}J_{Q|P}(\pi_v,\bs)$$
satisfy the properties  of Theorem 2.1 of \cite{A7}. 
If $v<\infty$ and $K_v$  is an open compact subgroup of $G(E_v)$, denote by 
$R_{Q|P}(\pi_v,\bs)_{K_v}$ the restriction of $R_{Q|P}(\pi_v,\bs)$ 
to the subspace $\H_P(\pi_v)^{K_v}$ of $K_v$-invariant vectors in the Hilbert 
space 
$\H_P(\pi_v)$ of the induced representation. If $v|\infty$, let $K_v\subset
G(E_v)$ be the standard maximal compact subgroup. For every $\sigma_v\in
\Pi(K_v)$ we denote by $\parallel\sigma_v\parallel$ the norm of the highest
weight of $\sigma_v$.  Given $\pi_v\in\Pi(G(E_v))$
and $\sigma_v\in\Pi(K_v)$, let $R_{Q|P}(\pi_v,\bs)_{\sigma_v}$ be the
restriction of $R_{Q|P}(\pi_v,\bs)$ to the $\sigma_v$-isotypical
subspace of $\H_P(\pi_v)$.  Finally for any place $v$, let $\Pi_{\di}(M(E_v))$
 be the subspace of all $\pi_v$ in
$\Pi(M(E_v))$ such that there exists an automorphic representation $\pi$ in the
discrete spectrum of $M(\A)$ whose local component at $v$ is equivalent to
 $\pi_v$. Then the main result that we need to prove Theorem \ref{th0.1} 
is the following proposition.

\begin{prop}\label{p0.2}
Let $v$ be a place of $E$. For all $M\in\cL$ and $P,Q\in\cP(M)$ the following 
holds.

1) If $v<\infty$, then for every  open compact subgroup $K_v$ 
of $\GL_n(E_v)$ and every multi-index $\alpha\in\N^r$ there exists $C>0$ such 
that
\begin{equation}\label{0.2}
\parallel D_{\bu}^\alpha R_{Q|P}(\pi_v,i\bu) _{K_v}\parallel \le C
\end{equation}
for all $\pi_v\in\Pi_{\di}(M(E_v))$ and $\bu\in\R^r$.

\smallskip
2) If $v|\infty$, then for every multi-index $\alpha\in\N^r$ there
 exist $C>0$ and $N\in\N$  such that
\begin{equation}\label{0.3}
\parallel D^\alpha_{\bu} R_{Q|P}(\pi_v,i\bu)_{\sigma_v}\parallel\le 
C(1+\parallel\sigma_v\parallel)^N
\end{equation}
for all $\bu\in\R^r$,  $\sigma_v\in\Pi(K_v)$ and 
 $\pi_v\in\Pi_{\di}(M(E_v))$.
\end{prop}

The normalization used in \cite{Mu4} differs slightly from the
normalization by $L$-functions. However,
it is easy to compare the two normalizations and it follows that 
Proposition \ref{p0.2} holds also with respect to the normalization used in
\cite{Mu4}. Together with Theorem 0.1 of \cite{Mu4}, this implies 
Theorem \ref{0.1}. Actually in \cite{Mu4} we considered only  reductive
 algebraic groups $G$ defined  over $\Q$. However, passing  to the group
$G^\prime=\Res_{E/\Q}G$ which is obtained from $G$ by restriction of scalars,
it follows immediately that the results of \cite{Mu4} can also be applied to
reductive algebraic groups defined over a number field.

The main analytic ingredients in the proof of Proposition
 \ref{p0.2} are a non-trivial uniform bound
 toward the Ramanujan 
hypothesis on the Langlands parameters of local components of cuspidal
automorphic representations \cite{LRS} and the determination of the residual 
spectrum \cite{MW}. 
Furthermore Corollary A.3  is important for the proof of (\ref{0.3}). 

Let us explain this in more detail.
First note that any local component $\pi_v$ of a cuspidal
automorphic representation $\pi$ of $\GL_m(\A)$ is generic \cite{Sk}. 
This implies that
$\pi_v$ is equivalent to a fully induced representation \cite{JS3}, i.e.,
$$\pi_v\cong I^{G(E_v)}_{P(E_v)}(\tau_1[t_1],...,\tau_r[t_r]),$$
where $P$ is a standard parabolic subgroup of type $(n_1,...,n_r)$, $\tau_i$
are  tempered representations of $\GL_{n_i}(E_v)$ and the $t_i$'s are 
real numbers satisfying
$$t_1>t_2>\cdots>t_r.$$
Here $\tau_i[t_i]$ is the representation 
$g\mapsto\tau_i(g)|\det(g)|^{t_i}$. For a unitary generic representation
$\pi_v$ the parameters $t_i$ satisfy $|t_i|<1/2$. In \cite{LRS}, Luo,
Rudnick and Sarnak proved that for an unramified $\pi_v$ which is the local 
component of a cuspidal automorphic representation of $\GL_m(\A)$, one has
\begin{equation}\label{0.5}
\max_i |t_i|<\frac{1}{2}-\frac{1}{m^2+1}.
\end{equation}
First we extend this result of Luo, Rudnick and Sarnak to all local 
components of cuspidal automorphic representations of $\GL_m(\A)$. Then we use
the description of the residual spectrum of $\GL_m(\A)$, given by M{\oe}glin
and Waldspurger \cite{MW}, to prove similar bounds  for the local
components of all automorphic representations in the residual spectrum of 
$\GL_m(\A)$ (cf. Proposition \ref{p3.4} for the precise statement). 
As a consequence, it follows that  for every local component $\pi_v$ of an 
automorphic
representation $\pi$ in the discrete spectrum of $M(\A)^1$ the normalized
intertwining operator $R_{Q|P}(\pi_v,\bs)$ is holomorphic in the domain
$\Re(s_i-s_j)> 2/(n^2+1)$, $1\le i<j\le r$. This is the key result which
is needed to prove Proposition \ref{p0.2}.
Combined with Corollary A.3  it immediately implies
(\ref{0.3}). For a finite place $v$ we use that by Theorem 2.1 of 
\cite{A7}  any matrix coefficient of
$R_{Q|P}(\pi_v,\bs)$ is a rational function of $q_v^{s_i-s_j}$, $i<j$. 
Together with the above result, this implies (\ref{0.2}). 

In an earlier version of this paper, the first two
authors were only able to establish (\ref{0.3}) for a fixed $K_v$-type, so
 that Theorem \ref{th0.1} could only be proved for $K$-finite functions
$f\in\Co^1(G(\A)^1)$. With the help of the appendix which was kindly provided
by E. Lapid, the $K$-finiteness assumption could be lifted. 
 
To extend the results of this paper to other reductive groups $G$ one would
need, in particular, the existence of non-trivial uniform bounds on the local
components of cuspidal automorphic representations of $G(\A)$. For a discussion
of this problem we refer to \cite{Sa}. Also note that \cite{CL} is a step in
 this direction.

The paper is organized as follows.
In section 2 we compare the two different normalizations of intertwining
operators and we prove some estimate for conductors. In section 3 we estimate
the (continuous) 
Langlands parameters of local components of cuspidal automorphic 
representations of $\GL_m$ which generalizes results of Luo, Rudnick and
Sarnak \cite{LRS} to the case of ramified representations. Then we use the
description of the residual spectrum of $\GL_m$ by M{\oe}glin and
Waldspurger \cite{MW} to obtain estimations for the Langlands parameters of
all local components of automorphic representations in the discrete spectrum
of $\GL_m$. We use these results in section 4 to prove Proposition \ref{p0.2}
and Theorem \ref{th0.1}. In the appendix, the normalized intertwining
operators for real Lie groups are studied. The main results is Corollary A.3 
which proves  estimations for  derivatives of matrix coefficients of
intertwining operators along
the imaginary axis, under the assumption that the intertwining operators are
holomorphic in a fixed strip containing the imaginary axis.

\smallskip
\noindent
{\bf Acknowledgment.} 
The first two  authors thank E. Lapid for many comments and suggestions which 
helped to
improve the paper considerably and for providing the appendix 
 by which the $K$-finiteness assumption in an earlier version of the 
paper could be lifted.

\section{Preliminaries}
\setcounter{equation}{0}

\subsection{}

Let $E$ be a number field and let $\A$ denote the ring of ad\`eles of $E$.
Fix a positive
integer $n$ and let $G$ be the group $\GL_n$ considered as
algebraic group over $E$. By a parabolic subgroup of $G$ we will
always mean a parabolic subgroup which is defined over $E$. Let
$P_0$ be the subgroup of upper triangular matrices of $G$. The Levi subgroup 
$M_0$ of $P_0$ is the group of diagonal matrices in $G$. A
parabolic subgroup $P$ of $G$ is called standard, if $P\supset
P_0$. By a Levi subgroup we will mean a subgroup of $G$ which contains $M_0$
and is the Levi component of a parabolic subgroup of $G$.
 If $M\subset L$ are Levi
subgroups, we denote the set of Levi subgroups of $L$ which contain $M$ by
${\cL}^L(M)$. Furthermore, let ${\cF}^L(M)$ denote the set of
parabolic subgroups of $L$ defined over $E$ which contain $M$, and let
${\cP}^L(M)$ be the set of groups in ${\cF}^L(M)$ for which $M$
is a Levi component. If $L=G$, we shall denote these sets by
${\cL}(M)$, ${\cF}(M)$ and ${\cP}(M)$. Write $\cL=\cL(M_0)$. Suppose that  
$P\in{\cF}^L(M)$. Then 
$$P=N_PM_P,$$
where $N_P$ is the unipotent radical of $P$ and $M_P$ is the unique Levi
component of $P$ which contains $M$. 

Suppose that $M\subset M_1\subset L$ are
Levi subgroups of $G$.  If $Q\in\cP^L(M_1)$ and $R\in{\cP}^{M_1}(M)$, there is a
unique group $Q(R)\in{\cP}^L(M)$ which is contained in $Q$ and whose
intersection with $M_1$ is $R$. 

Let $M\in\cL$ and denote by $A$  the split component of the center of $M$. 
Then $A$ is defined over $E$. Let $X(M)_E$ be the group of characters of
$M$ defined over $E$ and set 
\begin{equation*}
{\mathfrak a}_M=\mbox{Hom}(X(M)_E,{\mathbb R}).
\end{equation*}
Then $\af_M$ is a real vector space whose dimension equals that of $A$. Its
dual space is
\begin{equation*}
{\mathfrak a}^*_M= X(M)_E\otimes \R.
\end{equation*}

For any $M\in\L$ there exists a partition $(n_1,...,n_r)$ of $n$
such that $$M=\GL_{n_1}\times\cdots\times\GL_{n_r}.$$ Then $\af_M^*$
can be canonically identified with $(\R^r)^*$ and the Weyl group $W(\af_M)$
coincides with the group $S_r$ of permutations of the set
$\{1,...,r\}$.

\subsection{}
Let $H$ be a reductive algebraic group defined over $\Q$, let $F$ be a local  field of 
characteristic $0$ and let $K$ be an open compact subgroup of $H(F)$.
 We shall denote by $\Pi(H(\A))$ (resp. $\Pi(H(F))$, $\Pi(K)$, etc.) the set
of equivalence classes of irreducible unitary representations of $H(\A)$
(resp. $H(F)$, $K$, etc.).
\subsection{}

Let $F$ be a local field of characteristic zero. If $\pi$ is an admissible
representation of $\GL_m(F)$, we shall denote by $\widetilde\pi$ the
contragredient representation to $\pi$. 
Let $\pi_i$, $i=1,...,r$, be irreducible admissible
representations of the group $\GL_{n_i}(F)$. Then
$\pi=\pi_1\otimes\cdots\otimes\pi_r$ is an irreducible admissible
representation of
$$M(F)=\GL_{n_1}(F)\times\cdots\times\GL_{n_r}(F).$$ For
$\bs\in\C^r$ let $\pi_i[s_i]$ be the representation of
$\GL_{n_i}(F)$ which is defined by
$$\pi_i[s_i](g)=|\det(g)|^{s_i}\pi_i(g),\quad g\in\GL_{n_i}(F).$$
Let 
$$I_P^G(\pi,\bs)=
\Ind^{G(F)}_{P(F)}(\pi_1[s_1]\otimes\cdots\otimes\pi_r[s_r])$$ 
be the
induced representation and denote by $\H_P(\pi)$ the Hilbert space
of the representation $I_P^G(\pi,\bs)$. Sometimes we will denote 
$I^G_P(\pi,\bs)$ by $I^G_P(\pi_1[s_1],...,\pi_r[s_r])$.

\section{Normalizing factors for local intertwining operators}
\setcounter{equation}{0}

Let $F$ be a local field of characteristic 0. If $F$ is non-Archimedean, let
$\cO$ be the ring of integers of $F$ and let $\Pg$ be the unique maximal ideal
of $\cO$. Let $q$ be the number of elements of the residue field $\cO/\Pg$.
 Let $K=\GL_n(\cO)$. If $F$ is
Archimedean, let $K$ be the standard maximal compact subgroup of $\GL_n(F)$,
 i.e., $K=\rO(n)$, if $F=\R$, and $K=\U(n)$, if $F=\C$.

Let $M=\GL_{n_1}\times\cdots\times\GL_{n_r}$ be a standard Levi subgroup. We 
identify $\af_M$ with $\R^r$. Let $P_1,P_2\in\cP(M)$. Given $\pi\in\Pi(M(F))$,
 let
$$J_{P_2|P_1}(\pi,\bs),\quad\bs\in\C^r,$$
be the intertwining operator which intertwines
the induced representations $I_{P_1}^G(\pi,\bs)$ and
$I_{P_2}^G(\pi,\bs)$. The intertwining operator $J_{P_2|P_1}(\pi,\bs)$ is
defined by an integral over $N_{P_1}(F)\cap N_{\ov P_2}(F)$ which converges 
for $\Re(\bs)$ in a certain chamber of $\af_M^*$. 
It follows from \cite{A7} and \cite[\S15]{CLL} that the intertwining operators
 can be normalized in a suitable way. This means that there exist scalar 
valued meromorphic functions
$r_{P_2|P_1}(\pi,\bs)$ of $\bs\in\C^r$ such that the normalized
intertwining operators
$$R_{P_2|P_1}(\pi,\bs)=r_{P_2|P_1}(\pi,\bs)^{-1}J_{P_2|P_1}(\pi,\bs)$$
satisfy the properties of Theorem 2.1 of \cite{A7}.  The method used in
\cite{A7} works for every reductive group $G$. For $\GL_n$, however, it follows
from results of Shahidi \cite{Sh1}, \cite{Sh5} that local intertwining 
operators can be normalized by $L$-functions.

The normalizing factors defined by the Rankin-Selberg $L$-functions can be
 described as follows. Fix a nontrivial 
continuous character $\psi$ of the additive group $F^+$ of
$F$ and equip $F$ with the Haar measure which is selfdual with respect to 
$\psi$. First assume that $P$ is a standard maximal parabolic subgroup with 
Levi component $M=\GL_{n_1}\times\GL_{n_2}$. Let 
$\pi_i\in\Pi(\GL_{n_i}(F))$, $i=1,2$. 
If $F$ is non-Archimedean, let $L(s,\pi_1\times\pi_2)$ and
$\epsilon(s,\pi_1\times \pi_2,\psi)$ be the Rankin-Selberg $L$-function and
the $\epsilon$-factor, respectively,  as defined in \cite{JPS}.
If $F$ is Archimedean, let the $L$-function and the $\epsilon$-factor
be defined by using the Langlands parametrization (cf. \cite{A7}, \cite{Sh2}).
Then the normalizing factor can be regarded as a function 
$r_{\ov P|P}(\pi_1\otimes\pi_2,s)$ of one complex variable which is given by
\begin{equation}\label{2.1}
r_{\ov P|P}(\pi_1\otimes\pi_2,s)=\frac{L(s,\pi_1\times\widetilde\pi_2)}
{L(1+s,\pi_1\times\widetilde\pi_2)\epsilon(s,\pi_1\times\widetilde\pi_2,\psi)}.
\end{equation}
For arbitrary rank, the normalizing factors are products of normalizing factors
associated to rank one groups in $M$. 
Let  $e_i$, $i=1,\ldots,r$, denote the standard basis of $(\R^r)^*$. 
Then there exist
$\sigma_1,\sigma_2\in S_r$ such that the set of roots of $(P_1,A_M)$ and 
$(P_2,A_M)$, respectively, are given by
\begin{equation}\label{2.2}
\Sigma_{P_k}=\{ e_i-e_j\mid 1\le i,j\le r,\;\sigma_k(i)<\sigma_k(j)\},
\quad k=1,2.
\end{equation}
Put
\begin{equation*}
\I(\sigma_1,\sigma_2)=\{(i,j)\mid 1\le i,j\le r,\, \sigma_1(i)<
\sigma_1(j),\,\sigma_2(i)>\sigma_2(j)\}.
\end{equation*}
Then
$$\Sigma_{P_1}\cap\Sigma_{\ov P_2}=\{e_i-e_j\mid (i,j)\in
\I(\sigma_1,\sigma_2)\}.$$
Let  $\pi=\pi_1\otimes\cdots\otimes\pi_r$ where $\pi_i\in\Pi(\GL_{n_i}(F))$, 
$i=1,...,r$. For $\bs=(s_1,...,s_r)\in\C^r$ set
\begin{equation}\label{2.3}
\begin{split}
r&_{P_2|P_1}(\pi,\bs):  =\\
&\hskip10pt\prod_{(i,j)\in \I(\sigma_1,\sigma_2)}\frac{L(s_i-s_j,\pi_i\times
\widetilde{\pi}_j)}
{L(1+s_i-s_j,\pi_i\times \widetilde{\pi}_j)
\epsilon(s_i-s_j,\pi_i\times\widetilde{\pi}_j,\psi)}.
\end{split}
\end{equation}
Since the Rankin-Selberg $L$-factors are meromorphic functions, it follows that
$r_{P_2|P_1}(\pi,\bs)$ are meromorphic functions of $\bs\in\C^r$ and as 
explained in \cite[\S4]{A7} and \cite[p.87]{AC}, they satisfy all properties 
that are requested for normalizing factors.

In order to be able to apply the results of \cite{Mu4} we have to compare the
normalizing factors $r_{Q|P}(\pi,\bs)$ with those used in \cite{Mu4}
which we denote by $\widetilde r_{Q|P}(\pi,\bs)$. If $F$ is  Archimedean,
the normalizing factors $\widetilde r_{Q|P}(\pi,\bs)$ are defined
as the Artin $L$-factors and therefore, coincide with the
$r_{Q|P}(\pi,\bs)$. Assume that $F$ is non-Archimedean. By the
 construction of
the normalizing factors it suffices to consider the case where $P$ is
maximal, $Q=\ov P$ and $\pi$ is square integrable. Let $P$ be a standard
maximal parabolic subgroup of $\GL_m$ with Levi component
$$M=\GL_{m_1}\times\GL_{m_2}.$$
Then the normalizing factor
may be regarded as a function $\widetilde r_{\ov P|P}(\pi,s)$ of one complex
 variable $s$. 
We recall the construction of $\widetilde r_{\ov P|P}(\pi,s)$ for
square integrable representations $\pi$ \cite{CLL}. It follows from 
\cite{Si1}, \cite{Si2} that
for  every $\pi\in\Pi_2(M(F))$ there exists a rational function 
$U_P(\pi,z)$ such that the Plancherel measure
$\mu(\pi,s)$ is given by
$$\mu(\pi,s)=U_P(\pi,q^{-s}).$$
The rational function $U_P(\pi,z)$ is of the form
$$U_P(\pi,z)=a\prod_{i=1}^r\frac{(1-\alpha_iz)(1-\ov\alpha_i^{-1}z)}
{(1-\beta_iz)(1-\ov\beta_i^{-1}z)},$$
where $|\alpha_i|\le 1$, $|\beta_i|\le1$, $i=1,...,r$, and $a\in\C$ is a
constant such that
$$a\prod_{i=1}^r\frac{\alpha_i}{\beta_i}>0.$$
Let $b\in\C$ be such that
$$|b|^2\prod_{i=1}^r\frac{\ov\alpha_i}{\ov\beta_i}=a$$
and set
\begin{equation}\label{2.4}
V_P(\pi,z)=b\prod_{i=1}^r\frac{(1-\alpha_iz)}{(1-\beta_iz)}.
\end{equation}
Then the normalizing factor $\widetilde r_{\ov P|P}(\pi,s)$ is defined
by
$$\widetilde r_{\ov P|P}(\pi,s)=
V_P(\pi,q^{-s})^{-1}.$$
By definition we have
$$\mu_P(\pi,s)=\left(\ov{\widetilde r_{\ov P|P}(\pi,-\ov s)}\widetilde r_{\ov P|P}(\pi,s)\right)^{-1},$$
which is one of the main conditions that  the normalizing factors have to 
satisfy.

\smallskip
Let $\pi_1$ and $\pi_2$ be  tempered representations
 of $\GL_{m_1}(F)$ and $\GL_{m_2}(F)$, respectively.
By Corollary 6.1.2  of \cite{Sh1} the Plancherel measure is given by
$$\mu(\pi_1\otimes\pi_2,s)=q^{f(\widetilde{\pi}_1\times\pi_2)}
\frac{L(1+s,\pi_1\times\widetilde{\pi}_2)}{L(s,\pi_1\times\widetilde{\pi}_2)}
\frac{L(1-s,\widetilde{\pi}_1\times\pi_2)}{L(-s,\widetilde{\pi}_1\times\pi_2)},$$
where $f(\widetilde{\pi}_1\times\pi_2)\in\Z$ is the conductor of
$\widetilde{\pi}_1\times\pi_2$. Using the description of the Rankin-Selberg
$L$-functions for tempered representations \cite{JPS} (see also section 3),
it follows that 
\begin{equation}\label{2.5}
L(s,\pi_1\times\widetilde\pi_2)=\prod_{i=1}^{m_1}\prod_{j=1}^{m_2}
\left(1-a_{ij}q^{-s}\right)^{-1},
\end{equation}
with complex numbers $a_{ij}$ satisfying $|a_{ij}|<1$. Furthermore by
 \cite{JPS}, the
$\epsilon$-factor 
$\epsilon(s,\pi_1\times\pi_2,\psi)$ has the following form
\begin{equation}\label{2.6}
\epsilon(s,\pi_1\times\pi_2,\psi)=c(\pi_1\times\pi_2,\psi)q^{-f(\pi_1\times\pi_2,\psi)s},
\end{equation}

with $c(\pi_1\times\pi_2,\psi)\in\C-\{0\}$ and $f(\pi_1\times\pi_2,\psi)\in\Z$. Let
\begin{equation*}
c(\psi)=\max\{r\mid \Pg^{-r}\subset \ker\psi\}
\end{equation*}

Then
\begin{equation}\label{2.7}
f(\pi_1\times\pi_2,\psi)=n_1n_2c(\psi)+f(\pi_1\times\pi_2),
\end{equation}
with $f(\pi_1\times\pi_2)\in\Z$ independent of $\psi$. For simplicity assume
that $c(\psi)=0$. By Lemma 6.1 of \cite{Sh1} we have
\begin{equation}\label{2.8}
|c(\pi_1\times\pi_2)|=q^{f(\widetilde{\pi}_1\times\pi_2)/2}.
\end{equation}
Thus the $\epsilon$-factor can be written as
$$\epsilon(s,\pi_1\times\pi_2,\psi)=W(\pi_1\times\pi_2)
q^{(1/2-s)f(\pi_1\times\pi_2)},$$
where the root number $W(\pi_1\times\pi_2)$ satisfies $|W(\pi_1\times\pi_2)|
=1$. Finally, observe that $f(\pi_1\times\widetilde{\pi}_2)=f(\widetilde{\pi}_1\times
\pi_2)$. Using (\ref{2.5}), (\ref{2.6}) and (\ref{2.8}) it follows that
 the constant $b$ in (\ref{2.4}) can be chosen to be $\epsilon(0,\pi_1\times
\widetilde\pi_2,\psi)$ and
\begin{equation}\label{2.9}
\widetilde r_{\ov P|P}(\pi_1\otimes\pi_2,s)=
\frac{L(s,\pi_1\times\widetilde\pi_2)}{L(1+s,\pi_1\times\widetilde\pi_2)
\epsilon(0,\pi_1\times\widetilde\pi_2,\psi)}.
\end{equation}

Comparing (\ref{2.1}) and (\ref{2.9}), it follows that
\begin{equation}\label{2.10}
r_{\ov P|P}(\pi_1\otimes\pi_2,s)=
\frac{\epsilon(0,\pi_1\times\widetilde\pi_2,\psi)}
{\epsilon(s,\pi_1\times\widetilde\pi_2,\psi)}\widetilde r_{\ov P|P}(\pi_1\otimes\pi_2,s).
\end{equation}

This can be extended to parabolic groups of arbitrary rank in the usual way. 
Let $M$ be a standard Levi subgroup of $\GL_n$ of type $(n_1,...,n_r)$ and
let $P_1,P_2\in\cP(M)$. Using the 
product formula for $\widetilde r_{P_2|P_1}(\pi,\bs)$ \cite[p.29]{A7} and 
the corresponding product formula (\ref{2.3}) for $r_{P_2|P_1}(\pi,\bs)$, we
extend (\ref{2.10}) to all tempered representations $\pi$ of $M(F)$. 
Finally, if $\pi$ is any irreducible unitary representation of $M(F)$, it 
can be written as a Langlands quotient $\pi=J^M_R(\tau,\mu)$, where $R$ is
a parabolic subgroup of $M$, $\tau$ is a tempered representation of $M_R(F)$
and $\mu$ is a point in the chamber of $\af_R^*/\af_M^*$ attached to $R$. 
Then
$$r_{P_2|P_1}(\pi,\bs)=r_{P_2(R)|P_1(R)}(\tau,\bs+\mu)$$
and a similar formula holds for $\widetilde r_{P_2|P_1}(\pi,\bs)$
\cite[(2.3)]{A7}. Let  $\sigma_1,\sigma_2\in S_r$ be attached to $P_1,P_2$ 
such that $\Sigma_{P_i}$ is given by (\ref{2.2}).
Then we get
\begin{lem}\label{l2.1} 
For all irreducible unitary  representation $\pi=\otimes_{i=1}^r\pi_i$ of
$M(F)$ we have
\begin{equation*}
r_{P_2|P_1}(\pi,\bs)=\prod_{(i,j)\in I(\sigma_1,\sigma_2)}
\frac{\epsilon(0,\pi_i\times\widetilde\pi_j,\psi)}
{\epsilon(s_i-s_j,\pi_i\times\widetilde\pi_j,\psi)}
\widetilde r_{P_2|P_1}(\pi,\bs),\quad \bs\in\C^r.
\end{equation*}
\end{lem}

Since we will be concerned with logarithmic derivatives of normalizing factors,
we need estimates for $f(\pi_1\times\pi_2)$.
Let $f(\pi_i)$ be the conductor of $\pi_i$, $i=1,2$. Then by Theorem 1 of \cite{BH} and
Corollary (6.5) of \cite{BHK} we have
\begin{equation}\label{2.11}
0\le f(\pi_1\times\pi_2)\le n_1 f(\pi_1)+n_2 f(\pi_2)-\inf\{f(\pi_1),f(\pi_2)\}
\end{equation}

for all admissible smooth representations $\pi_i$ of $\GL_{n_i}(F)$, $i=1,2$.
Furthermore, by Corollary (6.5) of \cite{BHK}, we have
$f(\pi_1\times\pi_2)=0$ if and only if there exists a quasicharacter
$\chi$ of $F^\times$ such
that both $\pi_1\otimes\chi\circ\det$ and $\pi_2\otimes\chi^{-1}\circ\det$
are unramified principal series representations. 

By (\ref{2.11}) it suffices to estimate the conductors $f(\pi_i)$, $i=1,2$.

 Given an open compact subgroup
 $K\subset\GL_m(F)$, let
$$\Pi(\GL_m(F);K)=\big\{\pi\in\Pi(\GL_m(F))\;\big|\; \pi^K\not=\{0\}\big\}.$$

\begin{lem}\label{l2.2}
For every open compact subgroup $K$ of $\GL_m(F)$ there exists $C>0$
such that $f(\pi)\le C$ for all $\pi\in\Pi(\GL_m(F);K)$.
\end{lem}

\begin{proof} 
In the first step  we reduce the proof to the case of  square-integrable 
representations.
Let $\pi\in\Pi(\GL_m(F))$. Then there exist a parabolic subgroup $P$ of
$\GL_m$ of type $(m_1,...,m_r)$, tempered representations  $\tau_j$ of
$\GL_{m_j}(F)$ and real numbers $t_i$ with $t_1>\cdots>t_r$ 
such that $\pi$ is isomorphic to the Langlands quotient 
$J^{\GL_m}_P(\tau_1[t_1],...,\tau_r[t_r])$.  By Theorem 3.4 of \cite{J}
 it follows that
$$f(\pi)=f\left(I^{\GL_m}_P(\tau_1[t_1],...,\tau_r[t_r])\right)
=\sum_jf(\tau_j).$$
Furthermore a tempered representation $\tau$ of $\GL_d(F)$ is full 
induced: $\tau=I^{\GL_d}_Q(\sigma_1,...,\sigma_l)$,
where $Q$ is a parabolic subgroup of $\GL_d$ of type $(d_1,...,d_l)$ and
$\sigma_i$ is a square-integrable representation of $\GL_{d_i}(F)$, 
$i=1,...,l$. Then by (3.2.3) of \cite{J} we get
$$f(\tau)=\sum_jf(\sigma_j).$$
Next we relate the $K$-invariant subspaces.
We may assume that $K\subset\GL_m(\cO_F)$ is a congruence subgroup. 
Suppose that
$\pi$ is a subquotient of an induced representation $I^{\GL_m}_P(\sigma)$,
where $P$ is a parabolic subgroup of $\GL_m$ of type $(m_1,...,m_h)$, 
$\sigma=\otimes_i\rho_i$ and $\rho_i$ is an admissible
representation of $\GL_{m_i}(F)$.
If 
$\pi^K\not=\{0\}$, then $I^{\GL_m}_P(\sigma)^K\not=\{0\}$. Furthermore we have
\begin{equation*}
\begin{split}
I^{\GL_m}_P(\sigma)^K&=\left(I^{\GL_m(\cO_F)}_{\GL_m(\cO_F)\cap P}(\sigma)
\right)^K\\
&\hookrightarrow  \bigoplus_{\GL_m(\cO_F)/K}I^K_{K\cap P}(\sigma)^K\\
&\cong\bigoplus_{\GL_m(\cO_F)/K}\sigma^{K\cap P}.
\end{split}
\end{equation*}

Now observe that
$$K\cap P=\prod_{i=1}^hK_i,$$
where $K_i\subset\GL_{m_i}(\cO_F)$ are congruence subgroups. Then
$$\sigma^{K\cap P}\cong\otimes_{i=1}^h\rho_i^{K_i}.$$
Thus if $\pi^K\not=\{0\}$, then $\rho_i^{K_i}\not=\{0\}$ for all $i$, 
$1\le i\le h$. 
Combined  with  the above relations  of the conductors,
we reduce to the case of square-integrable representations. 

Let {\bf 1} denote the trivial representation of $K$.
By \cite[Theorem 10]{HC2} the set $\Pi_2(\GL_m(F),K)$ of 
square-integrable representations $\pi$ of $\GL_m(F)$ with 
$[\pi|_K:{\bf 1}]\ge 1$ 
is a compact subset of the space $\Pi_2(\GL_m(F))$ of square-integrable
representations of $\GL_m(F)$. By the definition of the topology in
$\Pi_2(\GL_m(F))$ \cite[\S2]{HC2}, the set
$\Pi_2(\GL_m(F),K)$ decomposes into a finite number of orbits 
under the canonical action of $i\R$ on $\Pi_2(\GL_m(F))$ given
by $\pi\mapsto\pi[it]$.
Since the conductor remains unchanged under twists by unramified characters,
the lemma follows.
\end{proof}

\section{Estimation of the Langlands parameters}
\setcounter{equation}{0}

For the unramified places, the Langlands parameters of local components of
cuspidal automorphic representations of $\GL_n(\A)$ have been estimated by
 Luo, Rudnick and Sarnak \cite{LRS}. The main purpose of this section is to 
extend the estimations of \cite{LRS} first to ramified places and then to 
local components of automorphic representations in the discrete spectrum of
$\GL_n(\A)$ in general. To deal with cuspidal automorphic representations
we follow the method of \cite{LRS} which uses properties of the Rankin-Selberg
$L$-functions.  First note that any local component of a cuspidal automorphic 
representation of $\GL_n$ is generic \cite{Sk}. Let $F$ be a local field. By
\cite{JS3}, any 
irreducible
generic representation $\pi$ of $\GL_n(F)$ is equivalent to  a fully 
induced representation
$$\pi= I^G_P(\sigma,\bs),$$
where $P$ is a standard parabolic subgroup of type $(n_1,...,n_r)$, 
$\bs=(s_1,...,s_r)\in\R^r$ satisfies $s_1\ge s_2\ge \cdots\ge s_r$ 
and $\sigma$ is a square integrable representation of $M_P(F)$. We shall refer
to $\bs$ as the (continuous) Langlands parameters of $\pi$. We also note that
 an irreducible induced representation $I^G_P(\sigma,\bs)$, $\sigma$ 
square-integrable and  $\bs\in\R^r$, is
unitary, only if it is equivalent to its hermitian dual representation 
$I^G_{\bar{P}}(\sigma,-\bs)$. By \cite[Theorem 7]{KZ} this implies that
there exists a $w \in W(a_M) $ of order 2 so that
$$(\sigma_1[s_1]\otimes\cdots\sigma_r[s_r])^w=\sigma_1[-s_1]\otimes
\cdots\otimes\sigma_r[-s_r].$$
Hence we have
\begin{equation}\label{3.1a}
\bigl\{\sigma_j[s_j]\bigr\}=\bigl\{\sigma_k[-s_k]\bigr\}.
\end{equation}
Moreover by the classification of the unitary dual of $\GL_n(F)$ it follows 
that
$$|\Re(s_i)|<1/2,\quad i=1,...,r.$$

The key result that we need about the local Rankin-Selberg $L$-functions
is the following lemma.

\begin{lem}\label{l3.1a}
Let $\pi$ be an irreducible unitary generic representation of $\GL_n(F)$, and
let $(s_1,...,s_r)\in\R^r$ be the Langlands parameters of $\pi$. Then 
$L(s,\pi\times\widetilde\pi)$ has a pole at the point
$$s_0=2\max_j|s_j|.$$
\end{lem} 

\begin{proof}
For the proof we need to describe the $L$-factors in more  detail.
Let $\pi\cong I^G_P(\sigma,\bs)$ be a fully induced representation with
$\sigma=\otimes_i\sigma_i$ for discrete series representations $\sigma_i$
of $\GL_{n_i}(F)$ and Langlands parameters $\bs=(s_1,...,s_r)$ satisfying 
the above conditions. Then by the multiplicativity of
the local $L$-factors \cite{Sh6} we get
\begin{equation}\label{3.1b}
L(s,\pi\times\widetilde\pi)=\prod_{i,j=1}^r L(s+s_i-s_j,
\sigma_i\times\widetilde\sigma_j).
\end{equation}
If $F$ is non-Archimedean, then this is Proposition 9.4 of \cite{JPS}. 
If $F$ is Archimedean, (\ref{3.1b}) follows from the Langlands classification
(see \S2 of \cite{Sh6}). 
This reduces the description of the $L$-factors to the case of
 square-integrable representations. We distinguish three cases according to 
the type of the field $F$.

\smallskip
{\bf 1.} $F=\R$.

The Rankin-Selberg local $L$-factors are defined in terms of $L$-factors
 attached to semisimple representations of the Weil group $W_\R$ by means 
of the Langlands correspondence \cite{L3}. If $\tau$ is a
semisimple representation of $W_\R$ of degree $n$ and $\pi(\tau)$ is the 
associated irreducible admissible representation of $\GL_n(\R)$, then
$$L(s,\pi(\tau))=L(s,\tau).$$
Furthermore, if
$$\tau=\bigoplus_{1\le j\le m}\tau_j$$
is the decomposition into irreducible representations of $W_\R$, then
$$L(s,\tau)=\prod_j L(s,\tau_j).$$
If $\tau'$ is another semisimple representation of $W_\R$ of degree $n'$ and
$\pi(\tau')$ is the associated irreducible admissible representation of
 $\GL_{n'}(\R)$, then the Rankin-Selberg local $L$-factor is given by
$$L(s,\pi(\tau)\times\pi(\tau'))=L(s,\tau\otimes\tau').$$
This reduces the computation of the $L$-factors to the case of irreducible
representations of the Weil group.

The irreducible representations of the Weil group $W_\R$ of  $\R$
are either 1 or 2 dimensional. The associated representations of $\GL_m(\R)$,
 $m=1,2$, are  square-integrable and all square-integrable representations
 are obtained in this way. Note that $\GL_m(\R)$
does not have square-integrable representations if $m\ge3$. 
To describe the $L$-factors, we define
Gamma factors by
$$\Gamma_\R(s)=\pi^{-s/2}\Gamma(s/2),\quad \Gamma_\C(s)=2(2\pi)^{-s}
\Gamma(s).$$
Suppose that $\tau$ is a two-dimensional
irreducible  representation of $W_\R$. Then
\[\tau=\mbox{ind} _{\C^*}^{W_\R}\theta,\]
where $\theta$ is a (not necessarily unitary)  character of $W_\C=\C^*$. Thus
there exist $t\in\C$ and $k\in\Z$ such that
\[\theta(z) = |z|^t(z/\ov z)^{k/2},\quad z\in\C^*.\] 
Then the $L$-factor is defined as
\[L(s,\tau ) = \Gamma_\C (s+t+ |k|/2).\]
The  one-dimensional irreducible representations  of $W_\R$ are of the form
\[\psi_{\epsilon, t} : (z,\sigma)\in W_\R\rightarrow \sign^{\epsilon}(\sigma)
|z|^t\] 
where $\epsilon =0,1$ and  ''sign'' is the sign
character of the Galois group.  The $L$-factor of $\psi_{\epsilon,t}$
is given by
\[L(s,\psi_{\epsilon,t}) =  
 \Gamma_\R(s+t+\epsilon).\]
Next we have to consider the tensor products of irreducible representations.

If $\tau =\mbox{ind} _{\C^*}^{W_\R}\theta$ and $\tau'= \mbox{ind}
_{\C^*}^{W_\R}\theta'$ are two-dimensional representations of $W_\R$,
 then we  have 
\begin{eqnarray*}
 \tau \otimes \tau'&=&
  \mbox{ind} _{\C^*}^{W_\R} (\theta \otimes
  \mbox{ind}_{\C^*}^{W_\R}{\theta'}|_{\C^*})\\
& =& \mbox{ind} _{\C^*}^{W_\R} (\theta \otimes {\theta'}) \oplus
\mbox{ind} _{\C^*}^{W_\R} (\theta \otimes {\theta'}^{\sigma}),
 \end{eqnarray*}
where $\sigma $ is the nontrivial element of the Galois group. Suppose that
$\theta'(z)=|z|^{t'}(z/\overline z)^{k'/2}$. Then we get
\begin{eqnarray*}\lefteqn{L(s,\tau\otimes\tau')}\\&
 =&
\Gamma_\C(s+t+t'+|k +k'|/2)\Gamma_\C(s+ t+t' +|k -k'|/2).\end{eqnarray*}
Similarly, if $\psi=\psi_{\epsilon,t'}$ is a one-dimensional representation,
then we have
$$\tau\otimes\psi=\mbox{ind}_{\C^*}^{W_\R}(\theta\otimes\psi|_{\C^*}),$$
and therefore, we get
$$L(s,\tau\otimes\psi)=\Gamma_\C(s+t+t'+|k|/2).$$
Finally, if $\psi=\psi_{\epsilon,t}$ and $\psi_{\epsilon',t'}$ are two
one-dimensional representations of $W_\R$, then
$$L(s,\psi\otimes\psi')=\Gamma_\R(s+t+t'+\widetilde\epsilon),$$
where $0\le\widetilde\epsilon\le 1$ and 
$\widetilde\epsilon=\epsilon+\epsilon'$ mod $2$. 

For $k\in\Z$ let  $D_k$ be the $k$-th discrete series representation of
$\GL_2(\R)$ with the same infinitesimal character as the $k$-dimensional
representation. Then $D_k$ is associated with the two-dimensional
representation  $\tau_k=\mbox{ind}^{W_\R}_{\C^*}(\theta_k)$ of $W_\R$ where the
character $\theta_k$ of $\C^*$ is defined by $\theta_k(z)=(z/\ov z)^{k/2}$.
For $\epsilon\in\{0,1\}$ let $\psi_\epsilon$ the character of
$\R^\times$, defined by $\psi_\epsilon(r)=(r/|r|)^\epsilon$. It corresponds to
the character $\psi_{\epsilon,0}$ of $W_\R$.
Using the above description of the $L$-factors, we get 
\begin{equation*}
\begin{aligned}
L(s,D_{k_1}\times D_{k_2})&=\Gamma_\C(s+|k_1-k_2|/2)\cdot
\Gamma_\C(s+|k_1+k_2|/2)\\
L(s,D_k\times\psi_\epsilon)&=L(s,\psi_\epsilon\times D_k)=\Gamma_\C(s+|k|/2)\\
L(s,\psi_{\epsilon_1}\times\psi_{\epsilon_2})&= \Gamma_\R(s+\epsilon_{1,2}), 
\end{aligned}
\end{equation*}
where $0 \leq \epsilon_{1,2}\leq 1$ with $ \epsilon_{1,2}\equiv
\epsilon _1 + \epsilon _2 \mbox{ mod }2$. Up to twists by unramified
characters this exhausts all possibilities for the $L$-factors in the 
square-integrable case.

\smallskip
{\bf 2.} $F=\C$.

As in the real case, the local $L$-factors are defined in terms of the 
$L$-factors attached to representations of the Weil group $W_\C$ by means of
the Langlands correspondence. The Weil group $W_\C$ is equal to $\C^*$.
Furthermore  we note that $\GL_m(\C)$ has square-integrable representations
 only if $m=1$. For $r\in\Z$ let $\chi_r$ be the character of $\C^*$ defined by
$\chi_r(z)=(z/\ov z)^r$, $z\in\C^*$. Then it follows that
$$L(s,\chi_{r_1}\times\chi_{r_2})=\Gamma_\C(s+|r_1+r_2|/2).$$
Again up to twists by unramified characters, these are  all possibilities
for the $L$-factors in the square-integrable case.

\bigskip
{\bf 3.} $F$ {\it non-Archimedean.}

Let $\pi$ be a square-integrable representation of
$\GL_m(F)$. 
By \cite{BZ} there is a divisor $d|m$,
a standard parabolic subgroup $P$ of $\GL_m(F)$ of type $(d,\ldots,d)$, and
an irreducible supercuspidal representation $\rho$ of $\GL_d(F)$ so that
$\pi$ is the unique quasi-square-integrable component of the induced
representation $I^G_P(\rho_1,\ldots,\rho_r)$, where $r=m/d$ and
$\rho_j=\rho\otimes |\det|^{j-(r+1)/2}$, $j=1,...,r$.
We will write $\pi=\Delta(r,\rho)$. The representation $\pi$
is unitary (or equivalently square-integrable) if only if $\rho$ is unitary.
Moreover the contragredient of $\Delta(r,\rho)$ is given by
$\widetilde{\Delta}(r,\rho)=\Delta(r,\widetilde{\rho})$.
Let $\sigma_1=\Delta(r_1,\rho_1)$ and $\sigma_2=\Delta(r_2,\rho_2)$
be square integrable representation of
$\GL_{m_1}(F)$ and $\GL_{m_2}(F)$, respectively. 
 Then by Theorem 8.2 of \cite{JPS} we have
$$L(s,\sigma_1\times\widetilde\sigma_2)=
\prod_{j=1}^{\min(m_1,m_2)}L(s+(m_1+m_2)/2-j,\rho_1\times\widetilde\rho_2).$$
Thus the description of the Rankin-Selberg $L$-functions is reduced to the
case of two supercuspidal representations.
Let $\rho_i$, $i=1,2$, be supercuspidal representations of $\GL_{k_i}(F)$. By
Proposition (8.1) of \cite{JPS} we have $L(s,\rho_1\times \rho_2)=1$ if
$k_1>k_2$. Since $L(s,\rho_1\times\rho_2)=L(s,\rho_2\times\rho_1)$, the same
holds for $k_1<k_2$. Let $k_1=k_2$. Then by Proposition (8.1) of \cite{JPS},
$L(s,\rho_1\times \rho_2)=1$, unless $\rho_1$ and $\rho_2$ are in the same
twist class, i.e., there exists $t\in\C$ such that $\rho_2\cong\rho_1[t]$. In
this case we have
\begin{equation*}
L(s,\rho_1\times \widetilde{\rho}_1[t])=L(s+t,\rho_1\times\widetilde{\rho}_1)=
\left(1-q^{-a(t+s)}\right)^{-1},
\end{equation*}
where $a|k_1$ is the order of the cyclic group of unramified characters
$\chi=|\det|^u$ such that $\rho_1\otimes\chi\cong\rho_1$. Let $\sigma=
\Delta(r,\rho)$. Then we get
\begin{equation*}
\begin{split}
L(s,\sigma\times\widetilde\sigma)&=\prod_{j=1}^{r}
L(s+r-j,\rho\times\widetilde\rho)\\
&=\prod_{j=1}^{r}\left(1-q^{-a(s+r-j)}\right)^{-1},
\end{split}
\end{equation*}
where $a$ is the order of the cyclic group of unramified characters
$\chi $ so that $\rho\otimes \chi $ is isomorphic to $\rho$.

From the above description of the local $L$-factors we conclude that they
have the following two properties. Let $\pi$ be a  square-integrable representation of $\GL_m(F)$.
Then  $L(s,\pi\times\widetilde\pi)$ has a pole at $s=0$. Furthermore, if 
$\pi_1$ and $\pi_2$ are square-integrable
representations of $\GL_{m_1}(F)$ and $\GL_{m_2}(F)$, respectively, then 
$L(s,\pi_1\times\pi_2)$ has no zeros. 

Now we are ready to prove the lemma. 
Let $\bs=(s_1,...,s_r)$ be the Langlands parameters of $\pi$. Let
$1\le i\le r$. Then it 
follows from (\ref{3.1a}) and (\ref{3.1b}) that 
$L(s,\pi\times\widetilde\pi)$ contains the factor $L(s-2s_i,\sigma_i\times
\widetilde\sigma_i)$. Using the above properties of the $L$-factors in the
square-integrable case, it follows that $L(s,\pi\times\widetilde\pi)$ has a
 pole at $2s_i$. By (\ref{3.1a}), $-s_i$ occurs also in $\bs$. Hence 
$L(s,\pi\times\widetilde\pi)$ has a pole at $2|s_i|$. In particular,
 $L(s,\pi\times\widetilde\pi)$ has a pole at $2\max_i|s_i|$. 
\end{proof}

Next we recall some facts about ray class characters. Let $E$ be a number
field. 
Let $\qf$ be a nonzero integral ideal of $E$ and denote by $C(\qf)$ the wide
ray class group of $E$ modulo $\qf$. We note that the term
''wide'' means that no positivity condition has been imposed at
the real places of $E$. Then a character of $C(\qf)$ is unramified
at all infinite places. Now recall that any character $\chi$ of
$C(\qf)$ can be identified with a character of the idele class
group $C_E=I_E/E^\times$ which is trivial on the congruence
subgroup $C_E^{\qf}=I_E^{\qf}E^\times/E^\times$ \cite{Neu}. Then
$\chi= \otimes_v\chi_v$ where $\chi_v$ is a character of
$E_v^\times$ which is unramified at all places $v|\infty$ and all
finite places $\pg\nmid\qf$. Furthermore for a finite place
$\pg\nmid\qf$, the character $\chi_{\pg}$ is given by
\begin{equation}\label{3.1}
\chi_{\pg}(\alpha)=\chi(\pg)^{v_{\pg}(\alpha)},\quad \alpha\in
E^\times_{\pg},
\end{equation}
where $v_{\pg}:E_{\pg}^\times\to\Z$ is the $\pg$-adic valuation.

Let $S$ be any finite set of finite  places of $E$. For an
integral ideal $\qf$ of $E$ let $X_{\qf}$ denote the set of all
wide ray class characters of conductor $\qf$ such that
$\chi(\pg)=1$ for all $\pg\in S$. Therefore we have $(\pg,\qf)=1$
for all $\pg\in S$. By (\ref{3.1}) it follows that
\begin{equation}\label{3.2}
\chi_{v}=1\quad\mbox{for all}\; v\in S\cup S_\infty.
\end{equation}
Let $X^*_{\qf}$ be the subset of $X_{\qf}$ consisting of all
primitive characters. Set $X=\cup_{\qf}X_{\qf}$ and
$X^*=\cup_{\qf}X^*_{\qf}$. For $\chi\in X^*_{\qf}$ and a cuspidal
automorphic representation $\pi$ of $\GL_n(\A)$, the partial
Rankin-Selberg $L$-function $L_S(s,(\pi\otimes\chi)\times
\widetilde\pi)$ is defined to be
$$L_S(s,(\pi\otimes\chi)\times\widetilde\pi)=\prod_{v\notin
S}L(s,(\pi_v\otimes\chi_v)\times\widetilde\pi_v).$$ 

We shall use the following  result of Luo, Rudnick and Sarnak which is
the main result of \cite{LRS}.

\begin{theo}\label{th3.1}
Given $n,\pi,S$ as the above and any $\beta>1-2/(n^2+1)$, there are
infinitely many $\chi\in X^*$ such that 
$$L_S(\beta,(\pi\otimes\chi)\times\widetilde\pi)\not=0.$$
\end{theo}

Now we can  establish our extension of Theorem 2 of \cite{LRS}.
\begin{prop} \label{p3.2}
Suppose that $\pi = \otimes_v \pi_v$ is a cuspidal automorphic 
representation of
$GL_n(\A)$ .  Let $s_v = (s_{1,v},\dots , s_{k,v})\in\R^k$ be
the Langlands parameters of the representation $\pi_{v}.$
  Then we have  
$$\max_j|s_{j,v}| <\frac{1}{2}-\frac{1}{n^2+1}.$$
\end{prop}

\begin{proof}
We follow the proof of Theorem 2 in \cite{LRS}.
Let $v$  be  a place of $E$  and set $S=\{v\}$.  Let $\chi\in X^*$, where
$X^*$ is the set of ray class characters with respect to $S$ which we defined 
above. The Rankin-Selberg $L$-function
 \[ L(s,(\pi \otimes \chi)\times\widetilde\pi) = L(s,(\pi_{v}\otimes\chi_v)
\times\widetilde\pi_{v}) L_S(s,(\pi\otimes\chi)\times\widetilde \pi).\]

is  holomorphic in the whole complex plane except for simple
poles at $s=1$ and $s=0$ if $\pi\otimes\chi\cong\pi$. 
 This follows from the work of
Jacquet, Piateski-Shapiro, Shalika, Shahidi, M{\oe}glin and Waldspurger
\cite{JPS}, \cite{JS1}, \cite{JS2}, \cite{MW}, \cite{Sh2}. Choosing the 
conductor of $\chi$ sufficiently large, we have $\pi\otimes\chi\not\cong\pi$.
Thus by
Theorem \ref{th3.1} we may choose $\chi\in X^*$ such that 
$L(s,(\pi\otimes\chi)\times\widetilde\pi)$ is an  entire function and by
(\ref{3.2}) we have $\chi_v=1$. 
Suppose that $s_0>0$  is a pole  of 
$$L(s,\pi_{v} \times \widetilde\pi_{v})=
L(s,(\pi_{v}\otimes\chi_v) \times \widetilde\pi_{v}).$$
Then $s_0$ must be a zero of 
$L_S(s,(\pi\otimes\chi)\times\widetilde \pi)$.  
Assume that $s_0>1- 2(1+n^2)^{-1}$. Then by  Theorem \ref{th3.1} 
there exists $\chi \in X^*$ with 
$L_S(s_0, (\pi\otimes\chi)\times\widetilde\pi)\not = 0$, $\chi_v=1$, and
$L(s, (\pi\otimes\chi)\times\widetilde\pi)$ entire.
Hence it follows that 
$s_0<1-2(1+n^2)^{-1}$. Together with Lemma \ref{l3.1a} the 
proposition follows.
\end{proof}

We shall now establish a similar result for the local components of residual 
automorphic representations of $\GL_n(\A)$. 
First we recall some facts about
representations of $\GL_n$ over a local field $F$. Any irreducible unitary
representation $\pi$ of $\GL_n(F)$ is equivalent to a Langlands quotient
$J^G_R(\tau,\mu)$. This is the unique irreducible quotient of an induced
representation $I^G_R(\tau,\mu)$ where $\tau$ is a tempered representation
of $M_R(F)$ and $\mu$ is a point in the positive chamber attached to $R$.
A slight variant of this description is as follows.
Recall that a tempered representation $\tau$ of $\GL_m(F)$ can be
described as follows. There exist a standard parabolic subgroup $Q$ of type 
$(m_1,...,m_p)$ and square-integrable representations $\delta_j$ of 
$\GL_{m_j}(F)$ so that
$\tau$ is isomorphic to the full induced representation
 $I^{\GL_m}_P(\delta_1,...,\delta_p)$. Hence by induction in stages,
there exist a standard parabolic subgroup $P$ of $G$ of type 
$(n_1,...,n_r)$,
discrete series representations $\delta_i$ of $\GL_{m_i}(F)$ and real
numbers $s_1\ge s_2\ge \cdots \ge s_r$ such that $\pi$ is equivalent to the 
unique irreducible quotient
$$J^{G}_P(\delta_1[s_1]\otimes\cdots\otimes\delta_r[s_r]),$$
 of the induced representation
$I^{G}_P(\delta_1[s_1]\otimes\cdots\otimes\delta_r[s_r])$ \cite[I.2]{MW}.
 We call $s_1,...,s_r$ the (continuous) Langlands parameters of $\pi$.

The residual spectrum for $\GL_n(\A)$ has been determined by M{\oe}glin and
Waldspurger \cite{MW}. Let $\pi=\otimes_v\pi_v$ be an irreducible automorphic
representation in the residual spectrum of $\GL_n(\A)$. By \cite{MW}, 
there is a divisor  $k|n$, a standard parabolic subgroup $P$ of type 
$(d,...,d)$, 
 and a cuspidal automorphic representation $\xi$ of $\GL_d(\A)$, $d=n/k$,
 so that the representation $\pi$ is a quotient of the induced representation
$$ I_{P(\A)}^{G(\A )}\left(\xi[(k-1)/2]\otimes
\cdots\otimes\xi[-(k-1)/2]\right).$$

\begin{lem}\label{l3.3}
 Let $\pi$ and $\xi$ be as above. Let $v$ be a place of $E$ and
let $s_1\ge\cdots\ge s_r$, be the
Langlands parameters of $\xi_v$. Then the Langlands parameters of $\pi_v$
are given by
$$\left(\frac{k-1}{2}+s_1,...,\frac{k-1}{2}+s_r,...,
-\frac{k-1}{2}+s_1,...,-\frac{k-1}{2}+s_r\right).$$
\end{lem}

\begin{proof}
Since $\xi_v$ is a local component of a cuspidal automorphic 
representation 
$\xi$ of $\GL_d(\A)$, it is generic. Using induction in stages, it follows
that  there exist a standard
parabolic subgroup $R$ of $\GL_d$ of type $(n_1,...,n_r)$,
discrete series representations $\delta_{i,v}$ of $\GL_{n_i}(E_v)$ and
real numbers $(s_1,...,s_r)$ satisfying
\begin{equation}\label{3.3}
s_1\ge s_2\ge \cdots\ge s_r;\quad|s_j|<1/2,\;
j=1,...,r.
\end{equation}
 such that $\xi_v $ is isomorphic to the full induced representation
$$\xi_v\cong I^{\GL_d}_R(\delta_{1,v}[s_1]\otimes\cdots\otimes
\delta_{r,v}[s_r]).$$

Let $Q =M_Q N_Q$
be the standard parabolic subgroup of $\GL_n$ whose Levi component $M_Q$
is a product of $k$ copies of $M_R$ and let 
$$\delta_v=(\delta_{1,v}\otimes\cdots\otimes\delta_{r,v})\otimes\cdots\otimes
(\delta_{1,v}\otimes\cdots\otimes\delta_{r,v}),$$
where $(\delta_{1,v}\otimes\cdots\otimes\delta_{r,v})$ occurs $k$ times.
 Define
\[\mu(k,\bs)=\left(\frac{k-1}{2} + s_1 ,\dots ,\frac{k-1}{2} +
s_r,\frac{k-3}{2} + s_1,\dots ,
 -\frac{k-1}{2}+s_r\right).\] 
By induction in stages we have
 \[I_P^G(\xi_v[(k-1)/2]\otimes\cdots\otimes
 \xi_v[-(k-1)/2]) = I^G_Q(\delta_v ,\mu(k,\bs)). \]
Furthermore, by (\ref{3.3})   the coordinates of $\mu(k,\bs)$ are
decreasing. Thus the induced representation 
$I^G_Q(\delta_v,\mu(k,\bs))$ has a
unique irreducible quotient which must be isomorphic to $\pi_v$.
\end{proof}

Next we recall a different method to parametrize irreducible unitary
 representations of $\GL_n(F)$. Let
$d|n$ and $k=n/d$. Let $P$ be the standard parabolic subgroup of type
$(d,...,d)$. Let $\delta$ be a discrete series representation of $\GL_d(F)$
and let $a,b\in\R$ be such that $b-a\in\N$. Then the induced representation
$$I^G_P(\delta[b]\otimes\delta[b-1]\otimes\cdots\otimes\delta[a])$$
has a unique irreducible quotient which we denote by $J(\delta,a,b)$.
 Especially, if $a=-(k-1)/2$ and $b=(k-1)/2$, then we put
\begin{equation}\label{3.4}
J(\delta,k):= J(\delta,a,b).
\end{equation}
By Theorem D of \cite{Ta} and \cite{Vo}, for every irreducible unitary 
representation of $\GL_n(F)$ there exist a standard parabolic subgroup $P$
of type $(n_1,...,n_r)$, $k_i|n_i$, discrete series representations 
$\delta_i$ of $\GL_{d_i}(F)$, $d_i=n_i/k_i$, and real numbers $s_1,...,s_r$
with $|s_i|<1/2$, $i=1,...,r$, such that $\pi$ is isomorphic to the fully 
induced representation:
\begin{equation}\label{3.5}
\pi\cong I^G_P(J(\delta_1,k_1)[s_1]\otimes\cdots\otimes J(\delta_r,k_r)[s_r]).
\end{equation}

Using this parametrization, we get the following analogue to Proposition \ref{p3.2} for local components of automorphic representations in the residual
 spectrum of $\GL_n(\A)$.

\begin{prop}\label{p3.4}
Let $\pi_v$ be a local component of an automorphic representation 
in the residual spectrum of $\GL_n(\A)$. There exist $k|n$,  a parabolic
 subgroup $P$ of type
 $(kn_1,\dots,kn_r)$, discrete series representations $\delta_{i,v}$ of 
$\GL_{n_i}(E_v)$ and real numbers $s_1,...,s_r$ satisfying
$$s_1\ge s_2\ge \cdots\ge s_r,\quad |s_i|<1/2 -(1+n^2)^{-1},\quad i=1,...,r,$$
such that 
$$\pi_v\cong I^G_P(J(\delta_{1,v}, k)[s_1]\otimes
\dots \otimes J(\delta_{r,v}, k)[s_r]).$$
\end{prop}

\begin{proof}
By the proof of Lemma \ref{l3.3}, $\pi_v$ is equivalent to a Langlands quotient of the
form
$$\pi_v\cong J^G_Q(\delta_v,\mu(k,\bs)),$$
where the parameters $\bs$ satisfy (\ref{3.3}).
Set
$$b_i=\frac{k-1}{2}+s_i,\qquad a_i=-\frac{k-1}{2}+s_i,\quad i=1,...,r.$$
Suppose there exist $1\le i<j\le r$ such that the triples 
$(\delta_{i,v},a_i,b_i)$ and $(\delta_{j,v},a_j,b_j)$ are linked in the sense
of I.6.3 or I.7 in \cite{MW}. Suppose that  $s_i\ge s_j$. Then it 
follows from (2) and (3)(i) on p.622  or from (1) and (2) on p.624 of \cite{MW}
that $a_i\ge a_j+1$ and $b_i\ge b_j+1$. This
implies $1\le |s_i-s_j|$ which contradicts (\ref{3.3}). 
Hence the triples $(\delta_{i,v},a_i,b_i)$ are pairwise not linked. Now observe
that 
$$J(\delta_{i,v},a_i,b_i)=J(\delta_{i,v},k)[s_i].$$
Let $P$ be the standard parabolic subgroup with Levi component
$$\GL_{kn_1}\times\cdots\times\GL_{kn_r}.$$
Let $\widetilde\delta_v=\otimes_{i=1}^r\delta_{i,v}$ and set
$$J(\widetilde\delta_v,k):=\otimes_{i=1}^rJ(\delta_{i,v},k).$$
Then it follows from Proposition I.9 of \cite{MW} that the induced 
representation
$$I^G_P(J(\widetilde\delta_v,k),\bs):=I^G_P(J(\delta_{1,v},k)[s_1]\otimes\cdots\otimes J(\delta_{r,v},k)[s_r])$$
is irreducible. By Lemma I.8, (ii), of \cite{MW}, 
$I^G_P(J(\widetilde\delta_v,k),\bs)$
is a quotient of the induced representation $I^G_Q(\delta_v,\mu(k,\bs))$. 
Since $I^G_P(J(\widetilde\delta_v,k),\bs)$ is irreducible, this quotient must
 be the
Langlands quotient. Thus
$$\pi_v\cong J^G_Q(\delta_v,\mu(k,\bs))\cong I^G_P(J(\delta_{1,v},k)[s_1]
\otimes\cdots\otimes J(\delta_{r,v},k)[s_r]).$$
By construction, the $s_i$'s are the Langlands parameters of a local component
of a cuspidal automorphic representation. Therefore it follows from Proposition
\ref{p3.2} that they satisfy
$$|s_i|<\frac{1}{2}-\frac{1}{n^2+1},\quad i=1,...,r.$$
\end{proof}

This result has an important consequence for the location of the poles of
normalized intertwining operators (see Proposition \ref{p4.2}).

\section{Proof of the main results}
\setcounter{equation}{0}

In this section we prove Proposition \ref{p0.2} and
Theorem \ref{th0.1}. 
To this end we need some preparation.

Let $M$ be a standard Levi subgroup  of type $(n_1,...,n_r)$  and let
 $P,P'\in\cP(M)$. Given  a place $v$ of $E$, let $\Pi_{\di}(M(E_v))$
denote the set of all $\pi_v\in\Pi(M(E_v))$ which are local components of
some automorphic representation $\pi$ in the discrete spectrum of $M(\A)$.
Without loss of generality, we may assume that $P$ is a standard parabolic 
subgroup.  
Let  $\pi_v\in\Pi_{\di}(M(E_v))$. Then $\pi_v=\otimes_i\pi_{i,v}$ with
  $\pi_{i,v}\in\Pi_{\di}(\GL_{n_i}(E_v))$, $i=1,...,r$. By
Proposition \ref{p3.4} there exist a standard parabolic subgroup $R_i$ of
$\GL_{n_i}$ with
$$M_{R_i}=\GL_{n_{i1}}\times\cdots\times\GL_{n_{im_i}},$$
$k_{ij}|n_{ij}$, discrete series representations $\delta_{ij}$ of 
$\GL_{d_{ij}}(E_v)$, $d_{ij}=n_{ij}/k_{ij}$, 
and $s_{ij}\in\R$ satisfying
\begin{equation}\label{4.8}
s_{i1}\ge\cdots\ge s_{im_i},\quad |s_{ij}|<
\frac{1}{2}-\frac{1}{n^2+1},
\end{equation}
such that
$$\pi_{i,v}\cong I^{\GL_{n_i}}_{R_i}\left(J(\delta_{i1},k_{i1})[s_{i1}]\otimes\cdots\otimes J(\delta_{im_i},k_{im_i})[s_{im_i}]\right).$$

Set $R=\prod_iR_i$. Then $R$ is a standard parabolic subgroup of $M$. Put
$m=m_1+\cdots+m_r$. We identify $\{(i,j)\mid i=1,...,r, \;j=1,...,m_i\}$ with
$\{1,...,m\}$ by
$$(i,j)\mapsto \sum_{k<i}m_k+j.$$
For $1\le l\le m$ let $(i,j)$ be the pair that corresponds to $l$. Put
$$\delta_l=\delta_{ij},\quad k_l=k_{ij},\quad s_l=s_{ij}.$$
Set
\begin{equation}\label{4.9}
J_{\pi_v}=\otimes_{i=1}^m J(\delta_i,k_i),\quad\bs_{\pi_v}=(s_1,...,s_m).
\end{equation}

Combing the above equivalences, we get
$$\pi_v\cong I^{M}_{R}(J_{\pi_v},\bs_{\pi_v}).$$

Let $P(R)$ and $P'(R)$ be the parabolic subgroups of $\GL_n$ with $P(R)\subset 
P$, $P'(R)\subset P'$, $P(R)\cap M=R$ and $P'(R)\cap M=R$. Then by induction in
stages, the induced representation $I^G_P(\pi_v,\bs)$ is unitarily equivalent
to the induced representation $I^G_{P(R)}(J_{\pi_v},\bs+\bs_{\pi_v})$. The
unitary map $\mu$ which provides this equivalence, is given by the evaluation
map as in \cite[p.31]{KS}.
Here we identify $\bs\in\C^r$ with an element in $\C^m$ by
$$(s_1,...,s_r)\mapsto (s_1,...,s_1,s_2,...,s_2,...,s_r,...,s_r),$$
where $s_i$ is repeated $m_i$ times.  \cite[p.31]{KS}.

\begin{lem}\label{l4.1}
Let $\mu$ be the unitary equivalence of the induced representations as above.
 Then we have
\begin{equation}\label{4.10}
\mu\circ R_{P'|P}(\pi_v,\bs)=R_{P'(R)|P(R)}(J_{\pi_v},\bs+\bs_{\pi_v})\circ\mu.
\end{equation}
\end{lem}

\begin{proof} 
First consider the unnormalized intertwining operators. Recall that the
unnormalized intertwining operators are defined by integrals which are
absolutely convergent in a certain shifted chamber \cite[Theorem 6.6]{KS}, 
\cite{Sh2}. If we compare
these integrals in their range of convergence as in \cite[p.31]{KS}, it follows
immediately that (\ref{4.10}) holds for the unnormalized intertwining operators.
So it remains to consider the normalizing factors. If we apply Lemma \ref{l3.3}
to each component $\pi_{i,v}$ of $\pi_v$, it follows that $\pi_v$ is equivalent
to a Langlands quotient of the form $J^M_Q(\delta_v,\mu(\bk,\bs_{\pi_v}))$,
where $\delta_v=\otimes_{i=1}^m\delta_i$ and $\mu(\bk,\bs_{\pi_v})=(\mu(k_1,\bs_{\pi_{1,v}}),...,\mu(k_r,\bs_{\pi_{r,v}})).$
Hence by (2.3) in \cite{A7} we get
$$r_{P'|P}(\pi_v,\bs)=r_{P'(Q)|P(Q)}(\delta_v,\bs+\mu(\bk,\bs_{\pi_v})).$$
Let $S_i$ be the standard parabolic subgroup of $M_{R_i}$ with Levi component
$$(\GL_{d_{i1}}\times\cdots\times\GL_{d_{i1}})\times\cdots\times(\GL_{d_{ip_i}}
\times\cdots\times\GL_{d_{ip_i}}),$$
where each factor $\GL_{d_{ij}}$ occurs $k_i$ times. Let $S=\prod_iS_i$. Then
$M_S$ and $M_Q$ are conjugate. Let
$w\in S_m$ be the element which conjugates $M_S$ into $M_Q$. Referring again 
to (2.3) in \cite{A7}, we get
$$r_{P'(R)|P(R)}(J_{\pi_v},\bs+\bs_{\pi_v})=r_{P'(R(S))|P(R(S))}(w\delta_v,\bs+
w\mu(\bk,\bs_{\pi_v})).$$
By (r.6) of \cite[p.172]{A8} it follows that
$$r_{P'(R)|P(R)}(J_{\pi_v},\bs+\bs_{\pi_v})
=r_{P'(R(S)^w)|P(R(S)^w)}(\delta_v,\bs+\mu(\bk,\bs_{\pi_v})),$$
where $R(S)^w=w^{-1}R(S)w$. Finally note that $P(Q)$, $P'(Q)$, $P(R(S)^w)$ and
$P'(R(S)^w)$ have the same Levi component and the reduced roots satisfy
$$\Sigma^r_{\ov{P'(R(S)^w)}}\cap\Sigma^r_{P(R(S)^w)}=
\Sigma^r_{\ov{P'(Q)}}\cap\Sigma^r_{P(Q)}.$$
Using the product formula (r.1) in \cite[p.171]{A8}, it follows that
$$r_{P'(R(S)^w)|P(R(S)^w)}(\delta_v,\bs+\mu(\bk,\bs_{\pi_v}))=
r_{P'(Q)|P(Q)}(\delta_v,\bs+\mu(\bk,\bs_{\pi_v})).$$
Combining the above equations, we get 
$$r_{P'|P}(\pi_v,\bs)=r_{P'(R)|P(R)}(J_{\pi_v},\bs+\bs_{\pi_v}),$$
and this finishes the proof of the lemma.
\end{proof}

We say that $ R_{Q|{P}}(\pi_v,\bs)$ has a pole at
$\bs_0\in\C^r$, if $ R_{Q|{P}}(\pi_v,\bs)$ has a matrix coefficient with
a pole at $\bs_0$. Otherwise, $R_{Q|{P}}(\pi_v,\bs)$ is called holomorphic in
$\bs_0$.

\begin{prop}\label{p4.2}  Let $M=\GL_{n_1}\times\cdots\times\GL_{n_r}$ be
 a standard Levi subgroup of $\GL_n$ and let $P,P^\prime\in\cP(M)$.
 For every place 
$v$ of $E$ and for all
$\pi_v\in\Pi_{\di}(M(E_v))$, the normalized intertwining operator
$ R_{P^\prime|{P}}(\pi_v,\bs)$ is holomorphic in the domain
$$\{\bs\in\C^r\mid \Re(s_i-s_j) >-2/(1+n^2),\; 1\le 1<j\le r\}.$$
\end{prop}

\begin{proof}
Using Lemma \ref{l4.1} we immediately reduce to the consideration of the
corresponding problem for
$R_{P^\prime(R)|P(R)}(J_{\pi_v},\bs+\bs_\pi)$, regarded as function of
 $\bs\in\C^r$. By Proposition I.10 of \cite{MW}, the intertwining operator
$R_{P^\prime(R)|P(R)}(J_{\pi_v},\bs)$
 is holomorphic in the domain of all $\bs\in\C^m$ satisfying
$\Re(s_i-s_j)>-1$, $1\le i<j\le m$.
Furthermore by (\ref{4.8}) the absolute value of all components  of
 $\bs_\pi$ is bounded by $1/2-(1+n^2)^{-1}$. Combining these observations
 the claimed result follows.

\end{proof}

\smallskip
\noindent
{\it Proof of Proposition \ref{p0.2}:}

Let $M$ be a standard Levi subgroup of type $(n_1,...,n_r)$ and let $P,Q\in
\cP(M)$. We distinguish between the Archimedean and non-Archimedean case. 

\smallskip
\noindent
{\bf Case 1:} $v|\infty$. 

Let $K_v$ be the standard maximal compact subgroup of $\GL_n(E_v)$. 
Given $\pi_v\in\Pi(M(E_v))$ and $\sigma_v\in\Pi(K_v)$, denote by 
$R_{Q|P}(\pi_v,\bs)_{\sigma_v}$ the restriction of $R_{Q|P}(\pi_v,\bs)$ to the 
$\sigma_v$-isotypical subspace $\H_P(\pi_v)_{\sigma_v}$ of the Hilbert space
$\H_P(\pi_v)$ of the induced representation. If $\pi_v$ belongs to 
$\Pi_{\di}(M(E_v))$, then by Proposition \ref{p4.2}, 
$R_{Q|P}(\pi_v,\bs)_{\sigma_v}$ has no poles in the domain of all
$\bs\in\C^r$ satisfying $|\Re(s_i)|<(1+n^2)^{-1}$, $i=1,...,r$. Using the 
factorization
of normalized intertwining operators and Corollary \ref{corrr} it follows that
there exist constants $C>0$ and $k\in\N$  such
that
\begin{equation*}
\parallel D_{\bu}R_{Q|P}(\pi_v,i\bu)_\sigma\parallel\le C
\left(1+\parallel \sigma_v\parallel\right)^k
\end{equation*}
for all $\pi_v\in\Pi(M(E_v))$, $\sigma_v\in\Pi(K_v)$ and $\bu\in\R^r$.
This proves (\ref{0.3}).

\smallskip
\noindent
{\bf Case 2:} $v<\infty$.

Using the properties of the normalized intertwining operators 
\cite[Theorem2.1]{A7}, one can factorize $R_{P^\prime|P}(\pi_v,\bs)$ in a 
product of normalized intertwining operators associated to maximal
parabolic subgroups. Thus we immediately reduce to the case where $P$ is
maximal and $P^\prime=\ov P$. Consider a matrix coefficient 
$(R_{\ov P|P}(\pi_v,s)v_1,v_2)$, where $\parallel v_1\parallel=\parallel v_2
\parallel=1$. By Theorem 2.1  of \cite{A7}, there is a rational
function $f(z)$ of one complex variable $z$ such that
\begin{equation}\label{4.4}
f(q^{-s})=(R_{\ov P|P}(\pi_v,s)v_1,v_2),\quad s\in\C.
\end{equation}
We shall now investigate the properties of the rational function $f$.
By Proposition I.10 of \cite{MW} we know that 
$(R_{\ov P|P}(\pi_v,s)v_1,v_2)$
is holomorphic in the half-plane $\Re(s)>0$. Hence $f(z)$ is holomorphic
in the punctured disc $0<|z|<1$. Moreover by unitarity of 
$R_{\ov P|P}(\pi_v,it)$, $t\in\R$, we have   
$|(R_{\ov P|P}(\pi_v,it)v_1,v_2)|\le 1$, $t\in\R,$
and hence $|f(z)|\le 1$ for $|z|=1$. To determine the behaviour of $f$ at
 $z=0$ we observe that the unnormalized intertwining operator 
$J_{\ov P|P}(\pi_v,s)$ is defined by an integral which is absolutely and
uniformly convergent in some half-plane $\Re(s)\ge c$.
Especially, $J_{\ov P|P}(\pi_v,s)$ is uniformly bounded for $\Re(s)\gg0$. 
The normalizing factor $r_{\ov P|P}(\pi_v,s)$ is given by (\ref{2.1}). It
follows from the expressions (\ref{2.5}) and (\ref{2.6}) for the 
$L$-factors and the epsilon factor, that there exist polynomials $P(z)$ and
$Q(z)$ with $P(0)=Q(0)=1$, a constant $a\in\C$ and $m\in\Z$ such that
$$r_{\ov P|P}(\pi_v,s)=aq^{ms}\frac{P(q^{-s})}{Q(q^{-s})},\quad s\in\C.$$
The integer $m$ is given by (\ref{2.7}) and it follows from (\ref{2.11}) that
there exists $c\ge0$, which depends on the choice of a nontrivial continuous
character of
$E^+_v$, such that $-c\le m$ for all $\pi_v\in\Pi(M(E_v))$.
Thus by the maximum principle it follows that for $0<|z|\le1$ we have
\begin{equation}\label{4.5}
|f(z)|\le \begin{cases}
                  1 & :\;m\ge0\\  |z|^m &:\; m<0.
                 \end{cases}
\end{equation}
Now assume that $\pi_v\in\Pi_{\di}(M(E_v))$. Then by Proposition \ref{p4.2},
$f(z)$ is actually holomorphic for $|z|<q^{2/(1+n^2)}$. Set
$$\delta=\min\{2,q^{2/(1+n^2)}\}.$$
Note that $\delta>1$. Let $\rho_1,...,\rho_r$ be the poles of $f$, where 
each pole is counted with its multiplicity. Let $-l$ be the order of $f$
at infinity. Set
$$g(z)=f(z)z^{-l}\prod_j\frac{z-\rho_j}{1-\ov\rho_jz}.$$
Since $|\rho_j|>\delta$, $j=1,...,r$, the rational function $g(z)$ is 
holomorphic for $|z|>1$, bounded for $|z|\ge1$ and satisfies
 $|g(z)|=|f(z)|\le 1$ for $|z|=1$. Thus $|g(z)|\le 1$ for $|z|\ge 1$ and
 hence,
$$|f(z)|\le |z|^l\prod_{j=1}^r\bigg|\frac{1-\ov\rho_jz}{z-\rho_j}\bigg|
=|z|^l\prod_{j=1}^r\bigg|\frac{\ov z-1/\rho_j}{1-z/\rho_j}\bigg|,\quad 
|z|\ge1.$$ 
For $1\le |z|<(1+\delta)/2$ the right hand side is bounded by 
$\left(\frac{1+\delta}{2}\right)^l\left(\frac{7}{\delta-1}\right)^r$. 
Together with (\ref{4.5}) it follows that there exists $C>0$, which is
independent of $\pi_v$, such that in the annulus 
$2/(1+\delta)< |z|<(1+\delta)/2$
we have
$$|f(z)|\le C\left(\frac{1+\delta}{2}\right)^l\left(\frac{7}{\delta-1}\right)^r.$$
Using Cauchy's formula we obtain a similar bound  for any derivative of $f$.
By (\ref{4.4}) this leads to a bound for any derivative of 
$(R_{\ov P|P}(\pi_v,s)v_1,v_2)$ in a strip $|\Re(s)|<\varepsilon$ for some
$\varepsilon>0$.
To complete the proof we need to  verify that for a given open compact
 subgroup $K_v$ of $\GL_n(E_v)$, the numbers $r$ and $l$ are bounded 
independently of $\pi_v$ if $\pi_v^{K_v\cap M(E_v)}\not=0$.

First consider $r$. By Theorem 2.2.2 of
\cite[p.323]{Sh2} there exists a polynomial $p(z)$ with $p(0)=1$ such that
$p(q^{-s})J_{\ov P|P}(\pi_v,s)$ is holomorphic on $\C$. Moreover the degree
of $p$ is bounded independently of $\pi_v$. Using the definition of the
 normalizing factors (\ref{2.1}), it follows immediately that there exists
a polynomial $\widetilde p(z)$ whose degree is bounded independently of
$\pi_v$ such that $\widetilde p(q^{-s})R_{\ov P|P}(\pi_v,s)$ is holomorphic
on $\C$. This proves that $r$ is bounded independently of $\pi_v$. 

To estimate $l$, we fix an open compact subgroup $K_v$ of $\GL_n(E_v)$. Our
goal is now to estimate the order at $\infty$ of any matrix coefficient of
$R_{\ov P|P}(\pi_v,s)_{K_v}$, regarded as a function of $z=q^{-s}$. Write 
$\pi_v$ as Langlands quotient $\pi_v=J^M_R(\delta_v,\mu)$ where $R$ is a 
parabolic subgroup of $M$, $\delta_v$ a square-integrable representation of
$M_R(E_v)$ and $\mu\in(\af_R^*/\af_M^*)_\C$ with $\Re(\mu)$ in the chamber 
attached to $R$.Then
$$R_{\ov P|P}(\pi_v,s)=R_{\ov P(R)|P(R)}(\delta_v,s+\mu)$$
with respect to the identifications described in \cite[p.30]{A7}. Here $s$ is
 identified with a point in $(\af_R^*/\af_G^*)_\C$ with repsect to the 
canonical embedding $\af_M^*\subset \af_R^*$. 
 Using again the factorization of normalized
intertwining operators we reduce to the case of a square-integrable
representation. 
Let ${\bf 1}$ denote the trivial representation of $K_v$.
By the same reasoning as in the proof of Lemma \ref{l2.2}
 we get
\begin{equation*}
[I^G_P(\delta_v,s)|_{K_v}:{\bf 1}]\le \#(\GL_n(\cO_v)/K_v)
[\delta_v|_{K_v\cap M(E_v)}:{\bf 1}].
\end{equation*}
By \cite[Theorem 10]{HC2} the set $\Pi_2(M(E_v),K_v)$ of 
square-integrable representations of $M(E_v)$ with 
$[\delta|_{K_v\cap M(E_v)}:{\bf 1}]\ge 1$ 
is a compact subset of the space $\Pi_2(M(E_v))$ of square-integrable
representations of $M(E_v)$. Under the canonical
action of $i\af_M^*$ in $\Pi_2(M(E_v))$, the subset
$\Pi_2(M(E_v),K_v)$ decomposes into a finite number of orbits.
In this way our problem is finally reduced to the consideration of the
matrix coefficients of $R_{\ov P|P}(\pi_v,s)_{K_v}$ for a finite number of
representations $\pi_v$. This implies the claimed bound for $l$.
\hfill$\Box$

\smallskip
\noindent
{\it Proof of Theorem \ref{th0.1}:}

Recall from \S2 that at finite places the normalization of the local 
intertwining operators differs from the normalization used in \cite{Mu4}.
Let $M=\GL_{n_1}\times\cdots\times\GL_{n_r}$, $Q,P\in\cP(M)$ and $v$
a finite place of $E$. Let $\widetilde r_{Q|P}(\pi_v,\bs)$, $\bs\in\C^r$,
be the normalizing factor used in \cite{Mu4} and let
$$\widetilde R_{Q|P}(\pi_v,\bs)=\widetilde r_{Q|P}(\pi_v,\bs)^{-1}
J_{Q|P}(\pi_v,\bs)$$
be the corresponding local normalized intertwining operator. Then it follows
from Lemma \ref{l2.1} together with (\ref{2.11}) and Lemma \ref{l2.2} that
for every multi-index $\alpha\in\N_0^r$ there exists $C>0$ such that
$$\parallel D^\alpha_{\bu}\widetilde R_{Q|P}(\pi_v,\bu)_{K_v}\parallel\le
C\sum_{|\beta|\le|\alpha|}\parallel D^\beta_{\bu} R_{Q|P}(\pi_v,\bu)_{K_v}
\parallel$$
for all $\bu\in\R^r$ and all $\pi_v\in\Pi_{\di}(M(E_v))$. Hence Proposition 
\ref{p0.2} holds also with respect to $\widetilde R_{Q|P}(\pi_v,\bs)$.
Together with Theorem 0.1 of \cite{Mu4} we obtain Theorem \ref{th0.1} of
the present paper.
\hfill$\Box$

\medskip
\noindent
{\bf Remark.} 
As the  proof of Proposition \ref{p0.2} shows,
the estimations (\ref{0.2}) and (\ref{0.3}) hold for all generic 
representations $\pi_v$ of $M(E_v)$
whose Langlands parameters $s_1,...,s_r$ satisfy a non-trivial bound of the
form $|s_i|<1/2-\varepsilon$, where $\varepsilon>0$ is independent of $\pi_v$. 
We note that this  assumption is really necessary and can not be removed in
 general. Especially,  
as the following example shows, the estimations can not be expected to 
be uniform in all $\pi_v\in\Pi(M(E_v))$.

\smallskip
\noindent
{\bf Example.}

Let $G=\GL_4(\R)$ and $P$ the standard parabolic subgroup with $M_P =
\GL_2 \times \GL_2 $. Consider the representation $I_P^G (\sigma
\times \sigma, \bs)$ where $\sigma$ is the spherical principal
series representation induced from the character $\mu = (\mu,
-\mu)$, $\mu$ real and $0 \leq  \mu \leq 1/2$ of the Borel
subgroup of $\GL_2(\R)$. We may assume that $\bs =(s,-s)$ with s
real. Then

\[I_P^G (\sigma \times \sigma, \bs) = I_B^G (\mu+s, -\mu +s,
\mu-s, \mu-s).\]

For fixed $0 \leq \mu <1/2 $ this representation is irreducible
for $0\leq |s| < 1/2-\mu $ and reducible for $|s| = 1/2 - \mu $
\cite{Sp}. The intertwining operator $R_{\bar{P}|P}(\sigma\times \sigma,
\bs)$ is therefore well defined on the interval $0\leq |s| <
1/2-\mu $ and has a pole for $-s = 1/2 - \mu $.
\hfill$\Box$

\smallskip
\noindent
The example shows that the poles of the normalized intertwining operator can
be arbitrary close to the imaginary axis. Thus, we can not expect to have
uniform bounds of the derivatives of the normalized intertwining operators
along the imaginary axis for all unitary $\pi$.

\begin{appendix}

\section{}

\setcounter{equation}{0}

\centerline{by Erez M. Lapid}

\newcommand{\bnd}{M}
Let $G$ be the real points of a connected reductive group defined
over $\mathbb{R}$. Let $K$ be a maximal compact subgroup of $G$
and let $P=M_PN_P$ be a parabolic subgroup of $G$ with its Levi
decomposition. Write $M=M_P=\,^0MA_M$ in the usual way. Let
$\sigma$ be an irreducible unitary representation of $\,^0M$
acting on a Hilbert space $H_\sigma$ and let $H_\sigma^\infty$ be
its smooth part. We denote by $I_\sigma^\infty$ the space of
smooth functions $f:K\rightarrow H_\sigma^\infty$ such that
$f(mk)=\sigma(m)f(k)$ for any $m\in K_M=M\cap K$ and $k\in K$ with
the inner product
\begin{equation*}
\langle f_1,f_2\rangle=\int_K(f_1(k),f_2(k))_{H_\sigma}\ dk
\end{equation*}
We denote the Lie algebra of $A_M$ by $\mathfrak{a}_M$. Let $P'$
be another parabolic subgroup of $G$ containing $M$ as its Levi
part. For any $\nu\in\mathfrak{a}_{M,\mathbb{C}}^*$, let
$J_{P'|P}(\nu)$ be the usual intertwining operator on
$I_\sigma^\infty$ (\cite[Chapter 10]{Wal2}) and let
$R_{P'|P}(\sigma,\nu)=r_{P'|P}(\sigma,\nu)^{-1}J_{P'|P}(\nu)$ be
the normalized intertwining operator (cf. \cite{A7}). Finally, for
any irreducible representation $\gamma$ of $K$ we denote by
$I_\sigma^G(\gamma)=I_\sigma(\gamma)$ the $\gamma$-isotypic part
of $I_\sigma^\infty$. We also denote by $\norm{\gamma}$ the norm
of the highest weight of $\gamma$.

The purpose of this appendix is to give a bound for the matrix
coefficients of the operator $R_{P'|P}(\sigma,\nu)$ on any
$K$-type near the unitary axis. By factoring
$R_{P'|P}(\sigma,\nu)$ it is enough to consider the ``basic'' case where
$P$, $P'$ are adjacent -- say along the root $\alpha$. In this
case the operator $J_{P'|P}(\sigma,\nu)$ depends only on
$(\nu,\alpha)$ and will be written as $J_{P'|P}(s)$ for
$(\nu,\alpha)=4s(\rho_P,\alpha)$. Similarly for
$R_{P'|P}(\sigma,s)$.

It follows from \cite[Lemma 10.1.11, Theorem 10.1.6,
10.1.13]{Wal2} that the poles of $J_{P'|P}(s)$ (counted with
multiplicities) are contained in
$\bigcup_{i=1}^r(\rho_i-\mathbb{N})$ for some complex numbers
$\rho_1,\dots,\rho_r$. By the nature of the normalization factors
we may enlarge the set $\{\rho_i\}$ to assume that the same holds
for $R_{P'|P}(\sigma,s)$ as well. Let
$$\bnd_\sigma^+=\max\{0,\Re\rho:\rho\text{ is a pole of }
R_{P'|P}(\sigma,s)\}
$$
and for any $\gamma\in\widehat{K}$ set
$$
\bnd_{\sigma,\gamma}^-=\max\{0,-\Re(\rho):\rho\text{ is a pole of
}R_{P'|P}(\sigma,s)\bigr\rvert_{I_\sigma(\gamma)}\}.$$

Finally, let
$$
\delta=\min\{\frac12,\abs{\Re(\rho)}:\rho\text{ is a pole of
}R_{P'|P}(\sigma,s)\}.
$$

\begin{lem} \label{main}

For any $\gamma\in\widehat{K}$ and any unit vectors
$\varphi_1,\varphi_2\in I_\sigma(\gamma)$ and any $\epsilon>0$ we
have
\[
\abs{(R_{P'|P}(\sigma,s)\varphi_1,\varphi_2)}\le
\left[(\bnd_\sigma^++\bnd_{\sigma,\gamma}^-+1)/\epsilon\right]^r
\]
in the strip $\abs{\Re(s)}<\delta-\epsilon$.
\end{lem}

\begin{proof}
Let $f(s)=(R_{P'|P}(\sigma,s)\varphi_1,\varphi_2)$. It is a
rational function of $s$ (\cite{A7}). We also have
$\abs{f(s)}\le1$ for $s\in i\mathbb{R}$ since $R_{P'|P}(\sigma,s)$
is unitary there. Define
\begin{equation*}
g(s)=f(s)\times\prod_{i=1}^r\prod_{j=\lceil\Re\rho_i+\delta\rceil}^{\lfloor
\Re\rho_i+\bnd_{\sigma,\gamma}^-\rfloor}
\frac{s-\rho_i+j}{s+\overline{\rho_i}-j}.
\end{equation*}
Then $g(s)$ is holomorphic (and rational) for $\Re(s)\le0$ and
$\abs{g(s)}=\abs{f(s)}\le1$ on $i\mathbb{R}$. Thus
$\abs{g(s)}\le1$ for $\Re(s)\le0$. It follows that in this region
\begin{equation*}
\abs{f(s)}\le\prod_{i=1}^r\prod_{j=\lceil\Re\rho_i+\delta\rceil}^{\lfloor
\Re\rho_i+\bnd_{\sigma,\gamma}^-\rfloor}
\abs{\frac{s+\overline{\rho_i}-j}{s-\rho_i+j}}\le
\prod_{i=1}^r\prod_{j=\lceil\Re\rho_i+\delta\rceil}^{\lfloor
\Re\rho_i+\bnd_{\sigma,\gamma}^-\rfloor}
\abs{\frac{\Re(s+\overline{\rho_i}-j)}{\Re(s-\rho_i+j)}}.
\end{equation*}
Then for $0\ge\Re s\ge-\delta+\epsilon$ each factor is bounded by
\begin{equation*}
\begin{split}
\prod_{j=\lceil\Re\rho_i+\delta\rceil}^{\lfloor\Re\rho_i
+\bnd_{\sigma,\gamma}^-\rfloor} \frac{-\Re s-\Re\rho_i+j}{\Re
s-\Re\rho_i+j}&<\prod_j\frac{\Re s-\Re\rho_i+j+1}{\Re
s-\Re\rho_i+j}\\
&<\frac{\bnd_{\sigma,\gamma}^-+1}{\epsilon}.
\end{split}
\end{equation*}
Similarly, one shows that for $0\le\Re s<\delta-\epsilon$
\begin{equation*}
\abs{f(s)}<\left[\frac{\bnd_\sigma^++1}{\epsilon}\right]^r.
\end{equation*}
\end{proof}



The following Proposition will be proved below.

\begin{prop} \label{nopole}
There exists a constant $c$ depending only on $G$ such that
\begin{equation} \label{maininequ}
\bnd_{\sigma}^+\le c, \  \ r\le c \, \ \ \bnd_{\sigma,\gamma}^-\le
c(1+\norm{\gamma})
\end{equation}
for all unitary $\sigma$.
\end{prop}

By Cauchy's formula, Lemma \ref{main} and Proposition \ref{nopole}
will imply the following.

\begin{corollary} \label{corrr}
For any differential operator $D(s)$ with constants coefficients
there exist constants $c'$, $k'$ (depending only on $G$) such that
\begin{equation}
 \label{bnd}
\norm{D(s)R_{P'|P}(\sigma,s)_{I_\sigma(\gamma)}}\le
c'\left(\frac{1+\norm{\gamma}}{\delta}\right)^{k'}
\end{equation}
for all $\gamma\in\widehat{K}$ and $s\in i\R$.
\end{corollary}

\begin{remark}
The example in \S 4 emphasizes that the dependence on $\delta$ is
essential if $\sigma$ is not tempered. This is already important
in order to lift the $K$-finiteness assumption in the absolute
convergence of the contribution of an individual cuspidal datum.
This point was overlooked in \cite{A4} (cf., p. 1329]). More
precisely, the property \cite[(7.6)]{A5} holds only for tempered
representations. We mention that it follows from \cite[Theorem
16.2]{KS} that for all $\sigma$ tempered we have $\delta>\delta_0$
where $\delta_0>0$ depends only on $G$.
\end{remark}

We will now prove Proposition \ref{nopole}. We first deal with the
first part of (\ref{maininequ}). More precisely, we have

\begin{lem} \label{s0}
There exists $s_0\in\R$, depending only on $G$, such that
$J_{P'|P}(\sigma,s)$ converges and $r_{P'|P}(\sigma,s)$ is
holomorphic and non-zero for all $\Re(s)>s_0$. In particular,
$\bnd_\sigma^+\le s_0$.
\end{lem}

\begin{proof}
Let $(Q,\tau,\lambda)$ be Langlands data for $\sigma$, i.e., $Q$
is a parabolic subgroup of $M$ with Levi subgroup $L$, $\tau$ is a
tempered representation of $L$ and $\lambda$ is a real parameter
in the positive Weyl chamber of $\mathfrak{a}_Q^*$ and $\sigma$ is
the irreducible quotient of the standard module defined by $Q$ and
$\lambda$. By \cite[5.5.2, 5.5.3]{Wal1}, or \cite[Ch. XI, Theorem
3.3]{BW} $\norm{\lambda}$ is bounded in terms of $G$ only.
Moreover identifying $I_\sigma^\infty$ with a quotient of
$I_\tau^\infty$, we may identify $J_{P'|P}(\sigma,s)$ with
$J_{QN_{P'}|QN_P}(\tau,s+\lambda)$ on the quotient space (cf.
\cite[p. 30]{A7} or \S 4). Moreover, we have
$r_{P'|P}(\sigma,s)=r_{QN_{P'}|QN_P}(\tau,s+\lambda)$. By
factoring $J_{QN_{P'}|QN_P}(\tau,s+\lambda)$ and
$r_{QN_{P'}|QN_P}(\tau,s+\lambda)$ the Lemma easily reduces to the
tempered case. Similarly, we reduce to the square-integrable case.
For $\sigma$ square-integrable we can take $s_0=0$ (\cite{A7}).
\end{proof}

The same argument reduces the second statement of
(\ref{maininequ}) to the square-integrable case. This case follows
from \cite[Theorem 16.2]{KS} and the compatibility of the
normalization factors with Artin's factors (\cite{A7}).

To continue the proof of Proposition \ref{nopole} we suppress for
the moment the assumption that $P$, $P'$ are adjacent and set
$\Sigma(P'|P)=\Sigma(P)\cap\Sigma(\overline{P'})$ where
$\overline{P'}$ is the parabolic opposite to $P'$ and
$\Sigma(P)=\Sigma(P,A_M)$ be the set of reduced roots of $A_M$ in
$P$.

The main assertion is the following.

\begin{lem} \label{nonzero}
There exists a constant $d$ (depending only on $G$) such that for
any $\gamma\in\widehat{K}$ $J_{P'|P}(\nu)$ is holomorphic and
injective on $I_\sigma(\gamma)$ in the domain
\[
\{\nu\in\mathfrak{a}_{M,\mathbb{C}}^*:
\Re(\nu,\alpha)>d(1+\norm{\gamma})\text{ for all }
\alpha\in\Sigma(P'|P)\}.
\]
\end{lem}

The last inequality of (\ref{maininequ}) then follows from Lemma
\ref{s0}, Lemma \ref{nonzero} and the relation
\[
R_{P|P'}(\sigma,-s)R_{P'|P}(\sigma,s)=\id.
\]

It remains to prove Lemma \ref{nonzero}. Clearly we may assume, by
passing to the derived group, that $G$ is semisimple. We first
need some more notation. Let $P_0=\,^0M_0A_0N_0$ be a minimal
parabolic subgroup of $G$, contained in $P$, so that $^0M_0$ is
compact. Let $\mathfrak{t}$ be a maximal abelian subalgebra of
$^0\mathfrak{m}$ and let
$\mathfrak{h}=\mathfrak{t}\oplus\mathfrak{a}_0$ (a direct sum with
respect to the Killing form). Then $\mathfrak{h}_{\C}$ is a Cartan
subalgebra of $\mathfrak{g}_{\C}$ and the real vector space
$\mathfrak{h}_R$ spanned by the co-roots is
$i\mathfrak{t}+\mathfrak{a}_0$ (\cite[2.2.5]{Wal1}). The Weyl
group $W=W(\mathfrak{g}_{\C},\mathfrak{h}_{\C})$ acts on
$\mathfrak{h}_R$ as well as on $\mathfrak{h}_{\C}^*$. We identify
the characters of the center of the universal enveloping algebra
of $\mathfrak{g}_{\C}$ as $W$-orbits of $\mathfrak{h}_{\C}^*$ via
the Harish-Chandra isomorphism. A similar discussion applies to
$M$. We denote by $\chi_\sigma$ the infinitesimal character of
$\sigma$.

For any (finite dimensional) irreducible representation $\sigma'$
of $^0M_0$ and $\mu\in\mathfrak{a}_{0,\C}^*$ we denote by
$\pi_{\sigma',\mu}=\pi_{\sigma',\mu}^G$ the corresponding
principal series representation on $G$. Its infinitesimal
character is (the $W$-orbit) of $\chi_{\sigma'}+\mu$ where
$\chi_{\sigma'}\in i\mathfrak{t}^*$ (the infinitesimal character
of $\sigma'$) is the translate of the highest weight of $\sigma'$
by the half-sum of positive roots in $\mathfrak{m}_0$.

\begin{lem}
There exists a constant $c$ depending only on $G$ such that any
unitary representation $\sigma$ of $M$ can be embedded
(infinitesimally) as a subrepresentation of a (non-unitary)
$\pi_{\sigma',\mu}^M$ and $\norm{\Re\mu}\le c(1+\norm{\gamma})$
whenever $\Hom_{K_M}(\gamma,\sigma)\ne0$.
\end{lem}

\begin{proof}
Suppose first that $\sigma$ is square-integrable. Using the
Casselman subrepresentation Theorem (e.g. \cite[Ch. 4]{Wal1} or
\cite[Theorem 8.37]{Kn}) we may embed $\sigma$ in some principal
series $\pi_{\sigma',\mu}^M$. By comparing infinitesimal
characters we infer that $\mu\in\mathfrak{a}_0^*$ and
\[
\norm{\chi_\sigma}^2=\norm{\mu}^2+\norm{\chi_{\sigma'}}^2\ge
\norm{\mu}^2.
\]
On the other hand by \cite[p.398]{Wa2}, (cf. \cite[p. 258]{Wal2})
the square of the norm of any $K$-type of $\sigma$ is bounded
below, up to a fixed additive constant, by $\norm{\chi_\sigma}^2$.
The lemma follows in this case.

To treat the general case we use the Langlands classification
Theorem to imbed $\sigma$ in $S(\tau,\lambda)$ where $Q$ is a
parabolic subgroup with Levi subgroup $L$, $\tau$ is a square-
integrable representation of $L$, $\lambda$ is in the closed
negative Weyl chamber of $\mathfrak{a}_Q^*$ and $S(\tau,\lambda)$
is the corresponding induced representation. As in the proof of
Lemma \ref{s0} we have $\norm{\lambda}<C$ independently of
$\sigma$. All $K$-types of $S(\tau,\lambda)$ (and hence, of
$\sigma$) contain a $K_L$-type of $\tau$ in their restriction to
$K_L$. Hence, by induction in stages, the Lemma reduces to the
square integrable case.
\end{proof}

We will now reduce Lemma \ref{nonzero} to the case where $P$ is a
minimal parabolic of $G$.

Imbed $\sigma$ in $\pi^M_{\sigma',\mu}$ as in the Lemma and
suppose that $\Re(\mu+\nu,\beta)>0$ for all
$\beta\in\Sigma(P_0N_{P'}|P_0N_P)$. Then
$J_{P'|P}(\pi_{\sigma',\mu},\nu)$ can be identified with
$J_{P_0N_{P'}|P_0N_P}(\sigma',\mu+\nu)$ and it is given by an
absolutely convergent integral. Its restriction to
$I_\sigma^\infty$ is $J_{P'|P}(\sigma,\nu)$. Thus, in that region
the injectivity of $J_{P'|P}(\sigma,\nu)$ on $I_\sigma(\gamma)$
follows from that of $J_{P_0N_{P'}|P_0N_P}(\sigma',\mu+\nu)$. We
note that the restriction to $A_M$ defines a bijection
$\alpha\leftrightarrow\alpha'$ between $\Sigma(P_0N_{P'}|P_0N_P)$
and $\Sigma(P'|P)$, and we have $(\nu,\alpha)=(\nu,\alpha')$. The
reduction follows.

By factoring $J_{P'|P}$ as a product of ``basic'' intertwining
operators we may also assume that $P'$ is adjacent to $P$. Let
$Q=LV$ be the parabolic subgroup generated by $P$ and $P'$. Then
$L$ has rank one and it follows from the argument of
\cite[10.4.5]{Wal2} that $J_{P'|P}(\sigma,\gamma)$ is injective on
$I_\sigma(\gamma)$ if and only if $J^L_{P'\cap L|P\cap
L}(\sigma,\nu^L)$ is injective on $I_\sigma^L(\gamma')$ for all
$\gamma'\in\widehat{K_L}$ which occur in the restriction of
$\gamma$. We observe that $\norm{\gamma'}\le\norm{\gamma}$ for
such $\gamma'$. Hence, we reduce to the case where $G$ is of rank
one, $P$ is minimal and $P'=\bar{P}$. Once again we can assume
that $G$ is semisimple as well. From now on we assume that this is
the case.

For $\Re(s)>0$ the representation $\pi_{\sigma,s\alpha}$ is of
finite length and its Langlands quotient is given by the image of
$J_{P'|P}(\sigma,s\alpha)$. Thus, $J_{P'|P}(\sigma,s\alpha)$ is
not injective on $I_\sigma(\gamma)$ if and only if $\gamma$ occurs
in one of the subquotients of $\pi_{\sigma,s\alpha}$ other than
the Langlands quotient. Assume that this is the case and let
$\pi'$ be any such subquotient. Then by \cite[Corollary
5.5.3]{Wal1} the Langlands parameter of $\pi'$ is smaller than
that of $\pi$. Thus, either $\pi'$ is square-integrable or $\pi'$
can be imbedded in $\pi_{\sigma',s'}$ with $0\le\Re(s')<\Re(s)$.
In the first case, the infinitesimal character of $\pi'$ is in
$\mathfrak{h}_R^*$, i.e., $s\in\R$, and by (\cite[p. 398]{Wa2})
\[
C+\norm{\gamma}^2\ge\norm{\chi_{\pi'}}^2=\norm{\chi_\sigma}^2+s^2\norm{\alpha}^2
\ge s^2\norm{\alpha}^2
\]
for a certain constant $C$. It follows that $s$ is bounded by a
constant multiple of $\norm{\gamma}$. In the second case, we have
\begin{equation} \label{infit}
\chi_{\sigma'}+s'\alpha=w(\chi_\sigma+s\alpha)
\end{equation}
for some $w\in W$. Write $w\alpha=\xi\alpha+\beta$ with $\xi\in\R$
and $\beta\in i\mathfrak{t}^*$. If $\beta=0$ then $\xi=\pm1$, $w$
stabilizes $i\mathfrak{t}^*$ and we obtain $s=\pm s'$ -- a
contradiction. Thus, $\beta\ne0$. Projecting (\ref{infit}) onto
$i\mathfrak{t}^*$ we obtain
\[
\chi_{\sigma'}=(w\chi_\sigma)_{i\mathfrak{t}^*}+s\beta.
\]
On the other hand, since $\gamma$ occurs $\pi_{\sigma,s\alpha}$,
$\sigma$ occurs in the restriction of $\gamma$ to $^0M=K_M$ and
hence $\norm{\sigma}\le\norm{\gamma}$. Similarly,
$\norm{\sigma'}\le\norm{\gamma}$. Once again, it follows that
$\abs{s}$ is bounded by a constant multiple of $1+\norm{\gamma}$.
This concludes the proof of Lemma \ref{nonzero}.


\end{appendix}


\begin{thebibliography}{MMM}
\bibitem[A1]{A1} J. Arthur, {\it Eisenstein series and the trace formula,}
  Proc. Sympos. Pure Math., {\bf 33}, Part I, Amer. Math. Soc., Providence,
  R.I. (1979), 253--274.
\bibitem[A2]{A2} J. Arthur, {\it A trace formula for reductive groups I: terms
  associated to classes in $G(\Q)$,} Duke. Math. J. {\bf 45} (1978), 911--952.
\bibitem[A3]{A3} J. Arthur, {\it On a family of distributions obtained from
  Eisenstein series I: Applications of the Paley-Wiener theorem},
  Amer. J. Math. {\bf 104} (1982), 1243--1288.
\bibitem[A4]{A4} J. Arthur, {\it On a family of distributions obtained from
  Eisenstein series II: Explicit formulas},
  Amer. J. Math. {\bf 104} (1982), 1289--1336.
\bibitem[A5]{A5} J. Arthur, {\it The trace formula in invariant form,} Annals
  of Math. {\bf 114} (1981), 1--74.
\bibitem[A6]{A6} J. Arthur, {\it The invariant trace formula. II. Global
  theory}, J. Amer. Math. Society {\bf 1} (1988), 501--554.
\bibitem[A7]{A7} J. Arthur, {\it Intertwining operators and residues. I.
Weighted characters}, J. Funct. Analysis {\bf 84} (1989), 19--84.
\bibitem[A8]{A8} J. Arthur, {\it On the Fourier transforms of weighted
  orbital integrals}, J. reine angew. Math. {\bf 452} (1994), 163--217.
\bibitem[AC]{AC} J. Arthur and L. Clozel, {\it Simple algebras, base change,
and the advanced theory of the trace formula}, Annals Math. Studies {\bf 120},
Princeton Univ. Press, Princeton, NJ, 1989.
\bibitem[Bu]{Bu} C.J. Bushnell, {\it Induced representations of locally
profinite groups}, Journal of Algebra {\bf 134} (1990), 104--114.
\bibitem[BH]{BH} C.J. Bushnell and G. Henniart, {\it An upper bound on
conductors for pairs}, J. Number Theory {\bf 65} (1997), 183--196.
\bibitem[BHK]{BHK} C.J. Bushnell, G.M. Henniart, and P.C. Kutzko, {\it Local
Rankin-Selberg convolutions for $\GL_n$: Explicit conductor formula}, J. Amer.
Math. Soc. {\bf 11} (1998), 703--730.
\bibitem[BW]{BW} A. Borel and N. Wallach, {\it Continuous cohomology, discrete
subgroups, and representations of reductive groups}, Annals of Math. Stud.
{\bf 94}, Princeton Univ. Press, 1980.
\bibitem[BZ]{BZ} I.N. Bernstein and A.V. Zelevinsky, {\it Induced
representations of reductive $p$-adic groups, I.} Ann. Scient. \'Ec. Norm.
Sup. {\bf 10} (1977), 441--472.
\bibitem[CLL]{CLL} L. Clozel, J.P. Labesse, and R. Langlands, {\it Morning
seminar on the trace formula}, mimeoraphed notes, IAS, Princeton, 1984.
\bibitem[CL]{CL} L. Clozel, {\it D\'emonstration de la conjecture} $\tau$. 
Invent. Math. {\bf 151} (2003), 297-328.
\bibitem[DH]{DH} A. Deitmar and W. Hoffmann, {\it On limit multiplicities for
spaces of automorphic forms}, Canad. J. Math. {\bf 51} (1999), 952-976.
\bibitem[GS]{GS} S. Gelbard and F. Shahidi, {\it Boundedness of automorphic
$L$-functions in vertical strips,} Preprint, 2000.
\bibitem[GJ]{GJ} R. Godement and H. Jacquet, {\it Zeta functions of simple
algebras}. Lecture Notes in Math. {\bf 260}, Springer, Berlin-Heidelberg-New
York, 1972.
\bibitem[HC1]{HC1} Harish-Chandra, {\it Automorphic forms on semi-simple Lie
groups},
Lecture Notes in Math. {\bf 62}, Springer-Verlag, Berlin-Heidelberg-New York,
1968.
\bibitem[HC2]{HC2} Harish-Chandra, {\it The Plancherel formula for reductive
p-adic groups}, In: Collected papers, IV, pp. 353-367, Springer-Verlag,
New York-Berlin-Heidelberg, 1984.
\bibitem[J]{J} H. Jacquet, {\it Principal $L$-functions of the linear group},
 Proc. Symp. Pure Math. {\bf 33} (1979), 63-86.
\bibitem[JPS]{JPS} H. Jacquet, I.I. Piatetski-Shapiro, and J.A. Shalika,
{\it Rankin-Selberg convolutions,} Amer. J. Math. {\bf 105} (1983), 367-464.
\bibitem[JS1]{JS1} H. Jacquet and J.A. Shalika, {\it Euler products and
 classification of automorphic representations I}, Amer. J. Math. {\bf 103}
 (1981), 499-558.
\bibitem[JS2]{JS2} H. Jacquet and J.A. Shalika, {\it Rankin-Selberg
convolutions: Archimedean theory}, in: Festschrift in Honor of I.I.
Piatetski-Shapiro on the occasion of his sixtieth birthday, Israel Math. Conf.
Proc. {\bf 2}, Weizmann, Jerusalem, 1990, 125-208.
\bibitem[JS3]{JS3} H. Jacquet, J.A.  Shalika, {\it The Whittaker models of 
induced representations}, Pacific J. Math. {\bf 109} (1983), 107--120. 
\bibitem[Kn]{Kn} A.W. Knapp, {\it Representation theory of semisimple groups},
Princeton University Press, Princeton, New Jersey 1986.
\bibitem[KS]{KS} A.W. Knapp and E.M. Stein, {\it Intertwining operators for
semisimple Lie groups, II}, Inventiones math. {\bf 60} (1980), 9-84.
\bibitem[KZ]{KZ} A.W. Knapp and G. Zuckerman, {\it Classification theorems for
representations of semisimple Lie groups}, in Lecture Notes Math. Vol. 
{\bf 587}, pp.138-159, Springer-Verlag, Berlin-Heidelberg-New York, 1977.
\bibitem[Lb]{Lb} J.-P. Labesse, {\it Noninvariant base change identities},
M\'em. Soc. Math. France (N.S.), {\bf 61}.
\bibitem[L1]{L1} R.P. Langlands, {\it{On the functional equations satisfied by
Eisenstein series,}} Lecture Notes in Math. {\bf{544}}, Springer--Verlag,
Berlin--Heidelberg--New York, 1976.
\bibitem[L2]{L2} R.P. Langlands, {\it Eisenstein series}, In: {\it Proc. Symp. Pure Math.} Vol. {\bf 9}, A. M. S., Providence, R.I., 1966.
\bibitem[L3]{L3} R.P. Langlands, {\it On the classification of irreducible
representations of real algebraic groups}, Mathematical Surveys and Monographs
{\bf 31}, Amer. Math. Soc. 1989.
\bibitem[LRS]{LRS} W. Luo, Z. Rudnick, and P. Sarnak, {\it On the generalized
Ramanujan conjecture for} $\GL(n)$. Proc. Symp. Pure Math. {\bf 66}, Vol. II, 
pp. 301--310, Amer. Math. Soc., Providence, Rhode Island, 1999. 
\bibitem[MW]{MW} C. Moeglin et J.-L. Waldspurger, {\it Le spectre r\'esiduel
de GL(n)}, Ann. scient. \'Ec. Norm. Sup., $4^e$ s\'erie, t. {\bf 22} (1989),
605--674.
\bibitem[Mu1]{Mu1} W. M\"uller, {\it The trace class conjecture in the theory
  of automorphic forms}, Annals of Math. {\bf 130} (1989), 473--529.
\bibitem[Mu2]{Mu2} W. M\"uller, {\it The trace class conjecture in the theory
  of automorphic forms. II}, Geom. Funct. Analysis {\bf 8} (1998), 315--355.
\bibitem[Mu3]{Mu3} W. M\"uller, {\it On the singularities of residual
intertwining operators}, Geom. funct. anal. {\bf 10} (2000), 1118-1170.
\bibitem[Mu4]{Mu4} W. M\"uller, {\it On the spectral side of the Arthur
trace formula}, Geom. funct. anal. {\bf 12} (2002), 669-722.
\bibitem[Neu]{Neu} J. Neukirch, {\it Algebraic number theory}, 
Springer-Verlag, Berlin-Heidel\-berg-New York, 1998.
\bibitem[RS]{RS} Z. Rudnick and P. Sarnak, {\it Zeros of principal 
$L$-functions and random matrix theory}.  Duke Math. J. {\bf 81} (1996),
269--322.
\bibitem [Sa]{Sa} P. Sarnak, {\it Notes on the generalized Ramanujan
 conjectures}, Lecture slides, Summer School on
Harmonic Analysis, The Trace Formula and Shimura Varieties, Fields Institute,
June 2003.
\bibitem [Se]{Se} A. Selberg, {\it Harmonic analysis}, in ''Collected Papers'',
Vol. I, Springer-Verlagf, Berlin-Heidelberg-New York (1989), 626--674.
\bibitem[Sh1]{Sh1} F. Shahidi, {\it Fourier transforms of intertwining
operators and Plancherel measures for} $\GL_n$, Amer. J. Math {\bf 106}
(1984), 67--111.
\bibitem[Sh2]{Sh2} F. Shahidi, {\it On certain $L$-functions}, Amer. J. Math.
{\bf 103} (1981), 297--356.
\bibitem[Sh3]{Sh3} F. Shahidi, {\it Local coefficients and normalization of
intertwining operators for $\GL_n$}, Comp. Math. {\bf 48} (1983), 271--295.
\bibitem[Sh4]{Sh4} F. Shahidi, {\it Local coefficients as Artin factors for
real groups,} Duke Math. J. {\bf 52} (1985), 973--1007.
\bibitem[Sh5]{Sh5} F. Shahidi, {\it A proof of Langlands' conjecture on
Plancherel measures; Complementary series for $p$-adic groups}, Annals Math.
{\bf 132} (1990), 273--330.
\bibitem[Sh6]{Sh6} F. Shahidi, {\it On multiplicativity of local factors},
Israel Math. conf. Proc. Vol. {\bf 3} (1990), 279--289.
\bibitem[Sk]{Sk} J.A. Shalika, {\it The multiplicity one theorem for $\GL_n$},
Annals Math. {\bf 100} (1974), 171-193.
\bibitem[Si1]{Si1} A. Silberger, {\it Special representations of reductive
$p$-adic groups are not integrable}, Annals Math. {\bf 111} (1980), 571--587.
\bibitem[Si2]{Si2} A. Silberger, {\it On Harish-Chandra's $\mu$-function for
$p$-adic groups}, Trans. Amer. Math. Soc. {\bf 260} (1980), 113--121.
\bibitem[Si3]{Si3} A. Silberger, {\it The Langlands quotient theorem for
$p$-adic groups}, Math. Ann. {\bf 236} (1978), 95--104.
\bibitem[Sp]{Sp} B. Speh , {\it The unitary dual of $\GL(3,\R)$ and $\GL(4, \R)$,} Math. Ann 258, (1981/82), 113-133.
\bibitem[Ta]{Ta} M. Tadi\'c, {\it Classification of unitary representations in irreducible representations of general linear group (non-Archimedean case)}. 
Ann. Sci. École Norm. Sup. (4) {\bf 19} (1986), 335--382. 
\bibitem[Vo]{Vo} D.B. Vogan, {\it The unitary dual of $\GL(n)$ over an
Archimedean field}, Inventiones math. {\bf 83} (1986), 449--505.
\bibitem[VW]{VW} D.B. Vogan, Jr. and N.R. Wallach, {\it Intertwining
operators for real reductive Lie groups}, Advances in Math. {\bf 82} (1990),
 203--243.
\bibitem[Wal1]{Wal1} N.R. Wallach, {\it Real reductive groups. I}, Pure and
Applied Mathematics vol. 132, Academic Press Inc., Boston, MA,
1988.
\bibitem[Wal2]{Wal2} N.R. Wallach, {\it Real reductive groups. II}, Pure and
Applied Mathematics vol. 132, Academic Press Inc., Boston, MA, 1992.
\bibitem[Wal3]{Wal3}
Nolan~R. Wallach.
\newblock Representations of reductive {L}ie groups.
\newblock In {\em Automorphic forms, representations and $L$-functions (Proc.
  Sympos. Pure Math., Oregon State Univ., Corvallis, Ore., 1977), Part 1},
  Proc. Sympos. Pure Math., XXXIII, pages 71--86. Amer. Math. Soc., Providence,
  R.I., 1979.
\bibitem[Wa1]{Wa1} G. Warner, {\it Harmonic analysis on semi-simple Lie groups
  I}, Springer-Verlag, Berlin-Heidelberg-New York, 1972.
\bibitem[Wa2]{Wa2} G. Warner, {\it Harmonic analysis on semi-simple Lie groups
  II}, Springer-Verlag, Berlin-Heidelberg-New York, 1972.
\bibitem[Wh]{Wh} E.T. Whitaker and G.N. Watson, {\it A course of modern
analysis}, Cambridge, University Press, 4th edition, 1940.
\bibitem[Z]{Z} A.V. Zelevinsky, {\it Induced representations of reductive
$p$-adic groups II. On irreducible reprresentations of} $\GL_n$, Ann. Scient.
\'Ec. Norm. Sup. {\bf 13} (1980), 165--210.
\end{thebibliography}
\end{document}